\documentclass[11pt,reqno]{amsart}

\numberwithin{equation}{section}

\begin{document}
\pagestyle{myheadings}
\markboth{Feigin, Jimbo, Loktev, Miwa and Mukhin}
{Bosonic formulas}

\title[Bosonic formulas]{Bosonic formulas for $(k,l)$--admissible partitions}
\author{B. Feigin, M. Jimbo,  S. Loktev, T. Miwa
 and E. Mukhin}
\address{BF: Landau institute for Theoretical Physics, Chernogolovka,
142432, Russia}\email{feigin@feigin.mccme.ru}  
\address{MJ: Graduate School of Mathematical Sciences, University of
Tokyo, Tokyo 153-8914, Japan}\email{jimbomic@ms.u-tokyo.ac.jp}
\address{SL: Institute for Theoretical and Experemental Physics and 
Independent University of Moscow}\email{loktev@mccme.ru}
\address{TM: Division of Mathematics, Graduate School of Science, 
Kyoto University, Kyoto 606-8502
Japan}\email{tetsuji@kusm.kyoto-u.ac.jp}
\address{EM: Dept. of Mathematics, University of 
California, Berkeley, CA 94720}\email{mukhin@math.berkeley.edu}

\date{\today}

\setcounter{footnote}{0}\renewcommand{\thefootnote}{\arabic{footnote}}

\begin{abstract}
Bosonic formulas for generating series of partitions with 
certain restrictions are obtained by solving a set of
linear matrix $q$-difference equations. Some particular cases
are related to combinatorial problems arising from solvable lattice models,
representation theory and conformal field theory.
\end{abstract}

\maketitle


\renewcommand{\theequation}{\thesection.\arabic{equation}}
\def\theenumi{\roman{enumi}}
\def\labelenumi{(\theenumi)}
\font\teneufm=eufm10
\font\seveneufm=eufm7
\font\fiveeufm=eufm5
\newfam\eufmfam
\textfont\eufmfam=\teneufm
\scriptfont\eufmfam=\seveneufm
\scriptscriptfont\eufmfam=\fiveeufm
\def\frak#1{{\fam\eufmfam\relax#1}}
\let\goth\frak

\newcommand{\sfrac}[2]{{\textstyle \frac{#1}{#2}}}
\newcommand{\g}{{\goth{g}}}
\newcommand{\slth}{\widehat{\mbox{\twelveeufm sl}}_2} 
\newcommand{\slN}{\mbox{\twelveeufm sl}_N} 
\font\twelveeufm=eufm10 scaled\magstep1
\font\fourteeneufm=eufm10 scaled\magstep2    
\font\seventeeneufm=eufm10 scaled\magstep3   
\font\twentyoneeufm=eufm10 scaled\magstep4   
\newcommand{\slNBig}{\mbox{\seventeeneufm sl}_N} 
\newcommand{\Z}{{\mathbb Z}} 
\newcommand{\C}{{\mathbb C}} 
\newcommand{\R}{{\mathbb R}} 
\newcommand{\bs}{\boldsymbol}
\def\b{{\bf b}}
\newcommand{\F}{{\mathcal F}}
\newcommand{\la}{\lambda}
\newcommand{\Ft}{\widetilde{\mathcal F}}
\renewcommand{\H}{{\mathcal H}}
\newcommand{\Cc}{{\mathcal C}}
\newcommand{\Oc}{{\mathcal O}}
\newcommand{\nn}{\nonumber}
\newcommand{\Ref}[1]{(\ref{#1})}
\newcommand{\bea}{\begin{eqnarray}}
\newcommand{\ena}{\end{eqnarray}}
\newcommand{\be}{\begin{eqnarray*}}
\newcommand{\en}{\end{eqnarray*}}
\newcommand{\lb}[1]{\label{#1}}

\newcommand{\id}{{\rm id}}
\newcommand{\bra}[1]{\langle #1 |}        
\newcommand{\ket}[1]{{| #1 \rangle}}      

\newcommand{\Remark}{\medskip \noindent {\it Remark.}\quad}
\newcommand{\example}{\medskip \noindent {\it Example.}\quad}
\newtheorem{thm}{Theorem}[section]
\newtheorem{cor}[thm]{Corollary}
\newtheorem{prop}[thm]{Proposition}
\newtheorem{lem}[thm]{Lemma}
\newtheorem{dfn}[thm]{Definition}


\setcounter{footnote}{0}
\renewcommand{\thefootnote}{\arabic{footnote})}
\renewcommand{\arraystretch}{1.2}

\setcounter{section}{0}
\setcounter{equation}{0}

\section{Introduction}
Let $k,l,n$ be positive integers. A partition of length $n$
\begin{equation}
\la=(\la_1\geq\la_2\geq\dots\geq \la_n\geq 0)=(0^{x_0},1^{x_1},2^{x_2},\dots)\notag
\end{equation}
is called $(k,l)$--admissible if 
\begin{equation}
\la_i\geq\la_{i+k}+l\qquad (i=1,\dots,n-k).\notag 
\end{equation}
Then a partition $\la$ is $(k,l)$--admissible if and only if the set of
integers $(x_0,x_1,\dots)$ satisfies
\begin{eqnarray}
0\leq x_i\leq k, \qquad x_i+\cdots+x_{i+l-1}\leq k.\label{KL}
\end{eqnarray}
for all $i\geq 0$.

The sequences satisfying (\ref{KL}) have been
discussed in the literature. For example, the case $(k,l)=(1,2)$
is related to the Rogers-Ramanujan identities and
is used in the computation of $1$-point correlation functions
of the hard hexagon model \cite{Baxter}. The case where $k$ is general
and $l=2$ appeared in the representation theory of the affine
Lie algebra $\widehat{sl}_2$ in \cite{FS}. More recently,
the case $k=2$ and $l\not\equiv1\bmod3$ appeared in the study
of vertex operator algebras related to the $(2,1)$ primary field 
in conformal field theory \cite{FJM}.

In this paper, we will not discuss such applications but treat
the problem purely combinatorially for general $k$ and $l$.

Given a positive integer $N$, consider $(k,l)$--admissible partitions
of length $n$ with $\la_1<N$. Let $d^{(N)}_{k,l}(a,n)$ be the number of such
partitions of weight $a$. In other words,
$d^{(N)}_{k,l}(a,n)$ is the number of sequences $(x_0,\dots,x_{N-1},0,\cdots)$,
satisfying \Ref{KL} such that
\begin{eqnarray*}
\sum_j jx_j=a,\qquad \sum_jx_j=n.
\end{eqnarray*}
In this paper we study the generating series of $(k,l)$--admissible
partitions:
\begin{equation}
\label{1DSUM}
\chi_{k,l}(q,z)=\displaystyle
\lim_{N\rightarrow\infty}\chi_{k,l}^{(N)}(q,z), \qquad 
\chi_{k,l}^{(N)}(q,z)=
\sum_{a,n}d^{(N)}_{k,l}(a,n)q^az^n.
\end{equation}

There are two different types of formulas which express the generating series
of the kind (\ref{1DSUM}) (one dimensional configuration sums): fermionic and
bosonic formulas (see, e.g., \cite{BMS}). In this paper, we present
a bosonic formula in the infinite volume limit $N\rightarrow\infty$.
The result is expressed as a finite sum of $(l-1)$--fold series. 
Each series is a sum of terms parameterized by a set of $(l-1)$ non-negative
integers; each term has the form
\[
\frac{(-1)^\alpha q^\beta z^\gamma}
{(q)_{t_2}\cdots(q)_{t_l}(q^{t_2+\cdots+t_l}z)_\infty},
\]
where $\alpha,\beta,\gamma,t_i$ are some 
integers described by explicit formulas.

It is worth mentioning the difference between our approach and that used,
e.g., in \cite{BMS}. Both approaches use finite polynomials
$\chi_{k,l}^{(N)}$ and recurrence relations which connect $N$-th and
$(N+1)$-st polynomials. 
In \cite{BMS} this recurrence is obtained by adding $x_N$ to the 
sequences $(x_0,\dots,x_{N-1})$ 
with fixed boundary conditions at the {\em right end}, whereas 
in our approach, the recurrence is obtained by adding $x_0$ to
sequences $(x_1,\dots,x_N)$ with fixed conditions at the {\em left end}. 

In the latter approach of the right end recurrence, one can handle
the generating series in a single variable $q$ by setting $z=1$.
In the former approach of the left end recurrence,
we have the shift of $z$ to $qz$ due to the shift of
weight. The advantage of the left end reccurence is that it makes
sense in the limit $N\to\infty$, where it becomes a $q$-difference
matrix equation.

To solve the $q$-difference equation we use the following method.
We consider the sum (\ref{1DSUM}) as a sum over the polytope bounded
by the conditions (\ref{KL}). 
In such a problem a classical idea 
is to express the result as a sum of contributions from the extremal points
of the polytope. 
In general, it is not easy to enumerate all the extremal points
of a polytope in large dimensions and compute the corresponding contributions. 
However, in our problem, the $q$-difference
equation reduces this problem to the one-dimensional case: we 
construct the extremal points and the corresponding contributions 
for the $(N+1)$-st polytope from those for the $N$-th.

The number of extremal points thus constructed grows exponentially in $N$.
Miraculously, most terms cancel in pairs, 
and the final formulas contain only polynomially many terms. 
We do not understand the reason why this cancellation takes place. 
Because of technical restrictions,  
in this paper we give 
the result only in the limit $N\rightarrow\infty$. 

A similar result by a similar method for the $\widehat{sl}_2$
coinvariants is obtained in \cite{JMM}.

The plan of the paper is as follows. In Section \ref{sec:2} we derive
the $q$-difference equations and define operators $A,B$ corresponding
to creation of new extremal points. In Section \ref{sec:3}
we determine the cancellation pairs. The proof of the cancellation is given
in Section \ref{sec:4}. Section 5 is devoted to describing the
remaining terms after the cancellation.

\setcounter{section}{1}
\setcounter{equation}{0}
\section{The $(k,l)$--configurations and their characters}
\lb{sec:2}

\subsection{The $(k,l)$--configurations}\lb{subsec:2.1}
Let $k,l$ be positive integers.
An infinite sequence of non-negative integers $\bs x
=(x_0,x_1,\dots)$ with finitely many non-zero entries
is called a {\em $(k,l)$--configuration} if 
for all $j\in\Z_{\ge 0}$ we have 
\begin{eqnarray}
x_j+\dots+x_{j+l-1}\le k.
\lb{config}
\end{eqnarray}
The smallest non-negative integer $N$ with 
the property $x_{i}=0$ for
all $i\ge N$ is called the {\em length} of the $(k,l)$--configuration
$\bs x$.

Let $\b=(b_0,\dots,b_{l-2})$ be a vector with non-negative integer
entries satisfying 
\begin{eqnarray}\label{b vector}
b_0\le b_1 \le\dots\le b_{l-2}\le k. 
\end{eqnarray}
Denote by $\mathcal{C}_{k,l;\b}^{(N)}$ the set of all
$(k,l)$--configurations of length at most $N$ such that 
$$
x_0\le b_0,~~x_0+x_1\le b_1,~\cdots,~
x_0+\dots+x_{l-2}\le b_{l-2}.
$$
We introduce the generating functions 
\begin{eqnarray}\label{(k,l)character}
\chi_{k,l;\b}^{(N)}(q,z)
=\sum_{\bs x\in\mathcal{C}_{k,l;\b}^{(N)}}
q^{\sum_{j\ge 0} jx_j}z^{\sum_{j\ge 0} x_j}, 
\end{eqnarray}
and call them the {\em characters of $(k,l)$--configurations}.

Setting $i=x_0$ and  
considering a shifted configuration $x'_j=x_{j+1}$,  
we are led to the following recursion equations.
\begin{lem}\label{recurrence lemma} 
The set of characters $\chi_{k,l;\b}^{(N)}(q,z)$ ($N=0,1,\dots$) 
is the unique solution of the recursion equations
\begin{eqnarray}
&&
\chi_{k,l;(b_0,\dots,b_{l-2})}^{(N+1)}(q,z)
=\sum_{i=0}^{b_0}z^i
\chi_{k,l;(b_1-i,b_2-i,\dots,b_{l-2}-i,k-i)}^{(N)}(q,qz),
\label{recurrence}
\end{eqnarray}
with the initial condition $\chi_{k,l;\b}^{(0)}(q,z)=1$ for all $b$.
\end{lem}

\subsection{The operators $A$, $B$}\lb{subsec:2.2}
Our main interest is to find a formula for the characters in ``bosonic
form''. 
To illustrate what it means, let us 
consider $N$ points $\{0,1,\cdots,N-1\}$ on a line segment $[0,N-1]$.  
Define its character to be the sum $1+z+\dots+z^{N-1}$. 
This polynomial can be written as a 
sum of two formal power series, 
\begin{eqnarray}\label{two end formula}
1+z+\dots+z^{N-1}=\frac1{1-z}+\frac{z^{N-1}}{1-z^{-1}}\,,
\end{eqnarray}
where the right hand side is understood as an expansion 
in non-negative powers of $z$. 
The two terms correspond to the two end points 
(``extremal points'') of the segment. 
The contribution {}from{} each extremal point may be viewed as 
the character of a half line (or a cone) $\{i\ge 0\}$, $\{i\le N-1\}$, 
respectively. 
We apply this principle to take the sum 
in the recursion equations \eqref{recurrence}.

Let $P_i$ ($i=1,\dots,l$) be monomials of $q,z$ 
of the form $q^{m_i}z^{n_i}$ where $m_i,n_{i}\in\Z_{\ge 0}$. 
Introduce the notation
$$
[P_1,\dots,P_l]= P_1^{b_0}P_2^{b_1-b_0}\dots
P_{l-1}^{b_{l-2}-b_{l-3}}P_l^{k-b_{l-2}}.
$$
We regard $[P_1,\dots,P_l]$ as a function of $q,z,\b,k$.
Note that 
$$
[\lambda P_1,\cdots,\lambda P_l]=\lambda^k[P_1,\cdots,P_l]. 
$$

Let $V_{l}$ be the complex vector space 
spanned by functions $\xi(q,z;\b,k)$ of the form $f(q,z)[P_1,\dots,P_l]$, 
where $f(q,z)$ is a formal power series in $q,z$ independent of $\b,k$, 
and $[P_1,\dots,P_l]$ belongs to the following list:
\begin{eqnarray*}
&&
q^{mk}z^{nk}[1,q^{c_1},q^{c_2},\cdots,q^{c_{l-1}}],
\\
&&
q^{mk}z^{nk}[q^{c_{l-1}}z,1,q^{c_1},\cdots,q^{c_{l-2}}],
\\
&&\cdots
\\
&&
q^{mk}z^{nk}[q^{c_{1}}z,q^{c_2}z,\cdots,q^{c_{l-1}}z,1],
\end{eqnarray*}
with  $c_i,m,n\in\Z_{\ge 0}$, $0\le c_1\le\dots\le c_{l-1}$.
We call such a function $\xi(q,z;\b,k)$ a simple vector,
$[P_1,\dots,P_l]$ the {\em vector part} and 
$f(q,z)$ the {\em scalar part} of it.

Let $S$ be the $q$-shift operator in $z$, 
defined by
\begin{eqnarray}
({S} g)(q,z)=g(q,qz).
\lb{shift}
\end{eqnarray}
Define two linear operators $A,B$ on $V_{l}$ by the formulas
\begin{eqnarray}
A(f(q,z)[P_1,\dots,P_l])&=&{{S}}\left(\frac{f(q,z)}{1-q^{-1}zP_l/P_1}\,
[P_1,P_1,P_2,\dots,P_{l-1}]\right),\label{OPA}
\\
B(f(q,z)[P_1,\dots,P_l])&=&{{S}}\left(\frac{f(q,z)}{1-qz^{-1}P_1/P_l}\,
[q^{-1}zP_l,P_1,P_2,\dots,P_{l-1}]\right).\label{OPB}
\end{eqnarray}
{}From{} the definition of the space $V_l$, it follows that  
${S} (q^{-1} z P_l/P_1)=z{{S}}(P_l/P_1)$ has 
the form $zq^n$ ($n\ge 0$) or $q^m$ ($m<0$). 
In the above, the factors 
${S}(1-q^{-1} z P_l/P_1)^{-1}$ 
and ${S}(1-q z^{-1} P_1/P_l)^{-1}$ 
are to be understood as power series expansions in $q,z$. 

It is immediate to check that
\begin{lem}\lb{lem:2.2}
The operators $A,B:\;V_{l}\to V_{l}$ are well defined.
\end{lem}
With the aid of $A,B$, 
the recursion equations \Ref{recurrence}  
can be recast into the following form. 

\begin{prop}\label{AB reccurence prop} 
For $N=0,1,\dots,$ we have 
$\chi_{k,l;\b}^{(N)}(q,z)\in V_{l}$. Moreover
$$
\chi_{k,l;\b}^{(N+1)}(q,z)=(A+B)\chi_{k,l;\b}^{(N)}(q,z).
$$
\end{prop}
\begin{proof}
Let $\xi(q,z;\b,k)=f(q,z)[P_1,\dots,P_l]$. 
Then using \Ref{two end formula}, we obtain:
\begin{align*}
& \sum_{i=0}^{b_0}z^i\xi(q,qz;b_1-i,\dots,b_{l-2}-i,k-i,k)
\\
& ={{S}} \left(f(q,z) P_1^{b_1}P_2^{b_2-b_1}\dots
P_{l-1}^{k-b_{l-2}}
\sum_{i=0}^{b_0}(q^{-1}z P_l/P_1)^i\right)
\\
&={{S}}\left(f(q,z) P_1^{b_1}P_2^{b_2-b_1}\dots
P_{l-1}^{k-b_{l-2}}
\left(\frac{1}{1-q^{-1}z P_l/P_1}+
\frac{(q^{-1}z P_l/P_1)^{b_0}}{1-qz^{-1}P_1/P_l}\right)\right)
\\
&=(A+B) \xi (q,z; \b,k).
\end{align*}
The assertion 
follows {}from{} linearity of operators $A,B$ 
and Lemma \ref{recurrence lemma}. 
\end{proof}
Note that for $N=0$ we can write
$$
\chi_{k,l;\b}^{(0)}(q,z)=1=[1,\dots,1]~~\in V_{l}.
$$
Let us denote the function $[1,\dots,1]\in V_l$ by $v_{\rm ini}$.
We then obtain the following expression for the characters. 
\begin{cor} The characters of $(k,l)$--configurations are given by
\begin{eqnarray}\label{trivial}
\chi_{k,l;\b}^{(N)}(q,z)=(A+B)^N v_{\rm ini}, \qquad N=0,1,\dots.
\end{eqnarray}
\end{cor}
\medskip

\noindent{\it Remark.}\quad
The operators $A$ and $B$ do not commute. 
It would be interesting to describe explicitly the algebra 
generated by $A,B$.

\subsection{Simplifications in the case $N=\infty$}\lb{subsec:2.3}
The tautological formula \Ref{trivial} 
gives the character as a sum over 
all monomials in $A,B$ of degree $N$. 
Here and below, by a {\em monomial} we mean an ordered composition of
operators of $A$, $B$, 
and we assign the degree to monomials by 
setting ${\rm deg}\,A={\rm deg}\,B=1$.
This expansion {\it a priori} comprises $2^N$ terms, 
so the number of terms grows exponentially with $N$. 
However we will see later (Section \ref{subsec:2.5}) that, 
at least in the limit $N\rightarrow\infty$, 
there are many cancellations and only polynomially many terms survive.   

{}From{} the definition of the characters \eqref{(k,l)character},
it is obvious that 
$$
\chi_{k,l;\b}^{(\infty)}(q,z)=
\lim_{N\rightarrow\infty}\chi^{(N)}_{k,l;\b}(q,z)
$$
exists as a formal power series in $q,z$. 
In the limit, the recursion equations \Ref{recurrence}  
turn into difference equations. 
\begin{prop}\label{infinite recurrence lemma}
The characters $\chi_{k,l;\b}^{(\infty)}(q,z)$ 
are the unique set of formal power series in $q,z$ 
satisfying the difference equation 
\begin{eqnarray}\label{infinite recurrence}
&&
\chi_{k,l;(b_0,\dots,b_{l-2})}^{(\infty)}(q,z)
=\sum_{i=0}^{b_0}z^i
\chi_{k,l;(b_1-i,b_2-i,\dots,b_{l-2}-i,k-i)}^{(\infty)}(q,qz),
\end{eqnarray}
and the initial condition 
$\chi_{k,l;\b}^{(\infty)}(q,0)=1$.
\end{prop}
\begin{proof}
We only need to show 
the uniqueness of the solution of \Ref{infinite
recurrence} with a given initial condition.

The difference equation \eqref{infinite recurrence} has the form 
\begin{eqnarray}\label{perturb rec}
y(z)=M(z) y(qz),
\end{eqnarray} 
where $y(z)$ is a vector whose components are labeled 
by $\b$ satisfying \Ref{b vector}, 
and $M(z)$ is a matrix whose non-zero entries are 
$z^i$ with some $i$. 
If we expand $y(z),M(z)$ as  
\begin{eqnarray*}
y(z)=\sum_{n\ge 0} y_n z^n,
\qquad
M(z)=\sum_{n\ge 0} M_n z^n,  
\end{eqnarray*}
then \eqref{perturb rec} reduces to the relations 
\begin{eqnarray*}
y_n=q^n M_0y_n+w_n
\qquad (n\ge 0), 
\end{eqnarray*}
where $w_n$ signifies terms which contain $y_m$ only with $m<n$. 
For $n>0$, $1-q^{n}M_0$ has an inverse matrix 
as a formal power series in $q$.  
Hence $y_n$ is uniquely determined by $y_0$. 
This proves the proposition.
\end{proof}

For the study of the difference equations 
\eqref{infinite recurrence}, 
it is helpful to first solve the much simpler equation 
\begin{eqnarray}
Av_\infty=v_\infty,
\qquad v_\infty\bigl|_{z=0}=v_{\rm ini}.
\lb{part}
\end{eqnarray}
Notice that for each $n$ we have
\begin{eqnarray}
A^n v_{\rm ini}=\frac{1}{(z)_n}
v_{\rm ini}, 
\lb{An}
\end{eqnarray}
where 
$(u)_n=\prod_{i=0}^{n-1}(1-q^iu)$.
Therefore the unique solution of \eqref{part} is given by 
\begin{eqnarray}\label{v infinity}
v_\infty=\lim_{N\rightarrow\infty}A^N v_{\rm ini}
=
\frac{1}{(z)_\infty} v_{\rm ini}=\frac{[1,\dots,1]}{(z)_\infty}\in V_l.
\end{eqnarray}

\begin{lem}\lb{lem:2.4}
We have 
\begin{eqnarray*}
&&\lim_{N\rightarrow\infty}(A+B)^{N}v_{\rm ini}
=v_\infty+\sum_M MBv_\infty
\\
&&=(1+B+AB+B^2+A^2B+AB^2+BAB+B^3+\dots)v_\infty,
\end{eqnarray*}
where the sum is over all monomials $M$ in $A,B$.
\end{lem}
\begin{proof}
The lemma will follow if we show that, 
for any sequence $C_1,\cdots,C_N$ of operators $A,B$, 
$\lim_{N\rightarrow\infty}C_1\cdots C_Nv_{\rm ini}=0$ 
holds unless $C_i=A$ for all sufficiently large $i$. 

Let us verify this statement. 
In view of \eqref{An}, 
it suffices to show that 
for any $C_i$ we have 
$\lim_{N\rightarrow\infty}C_1\cdots C_N B v_{\rm ini}=0$. 
The vector part of $Bv_{\rm ini}$ is $[z,1,\dots,1]$. 
Hence, for any monomial $M''$ of degree $l-1$, we have 
$M'' Bv_{\rm ini}=f[P_1,\cdots,P_l]$, 
where $f$ is a power series and each $P_i$ has degree $1$ in $z$. 
Applying further a monomial $M'$ of degree $d$, 
each component of the vector part of 
$M'M''Bv_{\rm ini}$ acquires a power $q^d$
due to the shift ${S}$. 
Therefore this expression vanishes when $d\rightarrow\infty$. 
\end{proof}

\subsection{Main result}\label{subsec:2.5}
Let $M,M'$ be monomials in $A,B$. 
We define an equivalence relation $M\sim M'$ if and only if 
there exist $r,r'\ge 0$ such that $MA^r=M'A^{r'}$. 
Lemma \ref{lem:2.4} is rephrased as 
\begin{eqnarray}
\chi^{(\infty)}_{k,l;\b}(q,z)=\sum_{M}Mv_\infty,
\lb{chkl}
\end{eqnarray}
where the sum ranges over all possible equivalence 
classes of monomials. 

We say that $M=C_1\cdots C_m$ 
is a {\em good monomial} if the following condition is 
satisfied.
\begin{eqnarray}
\mbox{For all $i$, $C_i=A$ implies $C_{i+l-1}=A$}.
\lb{goodM}
\end{eqnarray}
The notion of a good monomial 
makes sense also for its equivalence classes. 

\begin{example}
Here is the list of all good monomials 
up to equivalence for $l=2,3,4$. 
\begin{eqnarray*}
&&l=2~:~M=B^m\quad (m\ge 0),
\\
&&l=3~:~M=B^m (AB)^n \quad (m,n\ge 0),
\\
&&l=4~:~M=B^m (AB^2)^n (A^2B)^p\quad (m,n,p\ge 0),
\\
&&\phantom{l=4~:~M=}B^m (BAB)^n (A^2B)^p\quad 
(m\ge -1,n\ge 1,p\ge 0).
\end{eqnarray*} 
\end{example}

We are now in a position to formulate our main result.
\begin{thm}\lb{thm:2.1}
The character \eqref{chkl} is given by
\begin{eqnarray}
\chi^{(\infty)}_{k,l;\b}(q,z)=\sum_M{}^*
Mv_\infty,
\lb{chi*}
\end{eqnarray}
where the sum is taken over all 
the equivalence classes of good monomials. 
Moreover each term in the sum has the form 
\begin{eqnarray*}
&&\frac{(-1)^\alpha q^\beta z^\gamma}
{(q)_{t_2}\cdots(q)_{t_l}(q^{t_2+\cdots+t_l}z)_\infty},
\end{eqnarray*}
with some non-negative integers $\alpha,\beta,\gamma,t_2,\cdots,t_l$.
\end{thm}

Theorem \ref{thm:2.1} is proved in Sections \ref{sec:4}, \ref{sec:5}.

Compared to formula \Ref{chkl},  
Theorem \ref{thm:2.1} provides a big simplification 
because good monomials are much smaller in number among all monomials. 
In Section \ref{sec:5}, 
we give a parameterization of good
monomials for general $l$ and describe the
corresponding data $\alpha,\beta,\gamma,t_2,\cdots,t_l$ in combinatorial
terms. 
Below we write down the formula \eqref{chi*} explicitly in the cases
$l=2,3$.

\begin{example}
The case $l=2$:
\begin{eqnarray*}
\chi^{(\infty)}_{k,2;b_0}(q,z)
&=&\sum_{n\ge 0}(-1)^{n+1}q^{n^2k+nb_0+3n(n+1)/2}
z^{nk+b_0+n+1}
\frac{1}{(q)_n(q^nz)_\infty}
\\
&+&\sum_{n\ge 0}(-1)^{n}
q^{n^2k-nb_0+n(3n-1)/2}z^{nk+n}
\frac{1}{(q)_n(q^nz)_\infty}.
\end{eqnarray*}
This formula has been obtained in \cite{FL}.
\end{example}

\begin{example}
The case $l=3$:
Set 
\begin{eqnarray*}
g_{m,n}=q^{(3k+5)m^2/2+(3k+5)mn+2(k+2)n^2}z^{(k+1)(m+n)}.
\end{eqnarray*}
We retain the notation 
$[P_1,P_2,P_3]=P_1^{b_0}P_2^{b_1-b_0}P_3^{k-b_1}$. 
Then
\begin{eqnarray*}
\chi^{(\infty)}_{k,3;b_0,b_1}(q,z)=\sum_{m,n\ge 0}f_{m,n},
\end{eqnarray*}
where $f_{m,n}=B^m(AB)^nv_\infty$ is given as follows.
\begin{eqnarray*}
&&f_{3m,2n-1}=(-1)^m\frac{q^{-3(k+1)m/2-(k+2)n}}{(q)_m(q)_{m+2n-1}(q^{2m+2n-1}z)_\infty}
[1,q^{m},q^{-m-2n+1}z^{-1}]g_{m,n},
\\
&&
f_{3m,2n}
=(-1)^m\frac{q^{-(3k+1)m/2-kn}}{(q)_m(q)_{m+2n}(q^{2m+2n}z)_\infty}
[1,q^m,q^{2m+2n}]g_{m,n},
\\
&&
f_{3m+1,2n-1}
=(-1)^m\frac{q^{-(3k+1)m/2-(2k+1)n}}
{(q)_m(q)_{m+2n}(q^{2m+2n}z)_\infty}
[1,q^{m+2n},q^{2m+2n}]g_{m,n},
\\
&&
f_{3m+1,2n}
=(-1)^{m+1}
\frac{q^{(3k+5)m/2+(2k+3)n}z^{k+1}}
{(q)_m(q)_{m+2n}(q^{2m+2n}z)_\infty}
[1,q^{-2m-2n}z^{-1},q^{-m-2n}z^{-1}]g_{m,n},
\\
&&
f_{3m+2,2n-1}
=(-1)^{m+1}
\frac{q^{(3k+5)m/2+(k+2)n}z^{k+1}}
{(q)_m(q)_{m+2n}(q^{2m+2n}z)_\infty}
[1,q^{-2m-2n}z^{-1},q^{-m}z^{-1}]g_{m,n},
\\&&
f_{3m+2,2n}
=(-1)^{m+1}
\frac{q^{(3k+7)m/2+(k+4)n+1}z^{k+1}}
{(q)_m(q)_{m+2n+1}(q^{2m+2n+1}z)_\infty}
[1,q^{m+2n+1},q^{-m}z^{-1}]g_{m,n}.
\end{eqnarray*}
\end{example}

\subsection{Extremal configurations}\lb{subsec:2.4}
Before embarking upon the proof of Theorem \ref{thm:2.1}, 
let us make a remark about the ``extremal'' $(k,l)$--configuration 
associated with a given monomial in \Ref{chkl}.

Recall that the character \Ref{two end formula} of a ``segment''
is written as a sum of contributions {}from{} extremal points.
The same principle can be applied to write the character of a set of integer
points inside any convex ``polytope''in a higher dimensional space ${\bf R}^M$.
We expect that the result is a sum of rational functions over the vertices
(i.e., the extremal points) of the polytope, where each term represents 
the character of the cone corresponding to that vertex, see also
Lemmas 3.1, 3.2 in \cite{JMM}.
We will not discuss in how much generality this assertion is valid.
For a closely related discussion, see \cite{PKh}, Proposition 2.

Let us study this idea in our problem, where the set of integer points is
$\mathcal{C}_{k,l;\b}^{(N)}$ and $M=N$.
Take the simplest case $l=2$. We consider the generic situation $0<b_0<k$.

For $N=1$, the polytope is
a segment $\{x_0;0\leq x_0\leq b_0\}$.
We have $2$ vertices $a:=0$ and $b:=b_0$. They correspond to $Av_{\rm ini}$
and $Bv_{\rm ini}$, respectively.

For $N=2$, the polytope is a quadrilateral with the $4$ vertices
$aa:=(0,0),ab:=(0,k),ba:=(b_0,0),bb:=(b_0,k-b_0)$. They corresponds
to the terms $AAv_{\rm ini}$, etc.. We can classify these vertices in 
two groups, $aa$ and $ab$ with $x_0=0$ and $ba$ and $bb$ with $x_0=b_0$. 
These are the vertices of the sections of $\mathcal{C}^{(2)}_{b_0}$
by the hyperplanes $x_0=0$ and $x_0=b_0$. There exist natural bijections
between these sections and
$\mathcal{C}^{(1)}_k$ and $\mathcal{C}^{(1)}_{k-b_0}$, respectively.
The vertex $ab$ corresponds to a bad monomial $AB$. 
It turns out that in the limit $N\rightarrow\infty$ this vertex drops from 
the summation for the character formula. 
However, for $N=2$, we cannot drop it.

For $N=3$, the number of vertices is not equal to $8$ but is $7$:
This is because the hyperplane section $x_0=0$ is isomorphic to 
the degenerate polytope spanned by the integer points
$\mathcal{C}^{(2)}_0$, which is a triangle.
Both $aba$ and $abb$ correspond to a single vertex $(0,k,0)$.
The corresponding monomials $ABA$ and $ABB$ are bad, and 
two terms $ABAv_{\rm ini}$ and $ABBv_{\rm ini}$ cancel each other.
The degree $2$ bad monomial $AB$ finds a cancellation counterpart $ABB$
in the sense $ABv_\infty=ABAv_\infty=-ABBv_\infty$.

In Sections \ref{sec:3} and \ref{sec:4}, we will show
that such a cancellation always takes place in pairs. 

In the rest of this section we describe the extremal points
corresponding to good monomials. In Section \ref{sec:5}
the character formula is written as a sum over such vertices.

To a monomial $M=C_1\cdots C_N$, we associate a vertex
$(x_0,\dots,x_{N-1})\in\mathcal{C}_{k,l;\b}^{(N)}$ 
by the following inductive procedure. If $C_1=A$ we set $x_0=0$.
The section of $\mathcal{C}_{k,l;\b}^{(N)}$ by the hyperplane $x_0=0$
is isomorphic to $\mathcal{C}_{k,l;a(\b)}^{(N-1)}$ where
$a(\b)=(b_1,\ldots,b_{l-2},k)$. We set $(x_1,\ldots,x_{N-1})$ to be
the vertex corresponding to $C_2\cdots C_{N}$ in
$\mathcal{C}_{k,l;a(\b)}^{(N-1)}$. If $C_1=B$ we set $x_0=b_0$.
The section of $\mathcal{C}_{k,l;\b}^{(N)}$ by the hyperplane $x_0=b_0$
is isomorphic to $\mathcal{C}_{k,l;b(\b)}^{(N-1)}$ where
$b(\b)=(b_1-b_0,\ldots,b_{l-2}-b_0,k-b_0)$. We set $(x_1,\ldots,x_{N-1})$ to be
the vertex corresponding to $C_2\cdots C_{N-1}$ in
$\mathcal{C}_{k,l;b(\b)}^{(N-1)}$.
If a vertex corresponding to $M$ is $(x_0,\ldots,x_{N-1})$,
the vector part of the simple function
$M(v_{\rm ini})=C_1\cdots C_N(v_{\rm ini})$ is equal to
\begin{equation}\label{extr}
P_1^{b_0}P_2^{b_1-b_0}\cdots P_l^{k-b_{l-2}}=
\prod_{i=0}^{N-1}(q^iz)^{x_i}.
\end{equation}

We rephrase the above procedure as follows.

Consider an $l$-dimensional lattice $\Lambda$ generated by formal symbols 
$b_0,\ldots,b_{l-2},k$. We consider ${\bf Z}$ linear mappings $M_a,M_b$
on $\Lambda$ given by
\begin{eqnarray*}
M_a(b_i)&=&
\begin{cases}
b_{i+1}&\hbox{\rm if }0\leq i\leq l-3;\\
k&\hbox{\rm if }i=l-2,\\
\end{cases}\\
M_a(k)&=&k,\\
M_b(b_i)&=&
\begin{cases}
b_{i+1}-b_0&\hbox{\rm if }0\leq i\leq l-3;\\
k-b_0&\hbox{\rm if }i=l-2,
\end{cases}\\
M_b(k)&=&k.
\end{eqnarray*}

For a given monomial $M=C_1\cdots C_N$ we set
\[
c_i=
\begin{cases}
a&\hbox{\rm if }C_i=A;\\
b&\hbox{\rm if }C_i=B.
\end{cases}
\]

Then we have
\begin{equation}
x_i=
\begin{cases}
0&\hbox{\rm if } C_{i+1}=A;\\
M_{c_1}\circ\cdots\circ M_{c_i}(b_0)&\hbox{\rm if }C_{i+1}=B.
\end{cases}\label{XI}
\end{equation}

If $M$ is a good monomial, the above procedure is simplified. First we prepare
some notation and a few lemmas.

Define $\iota_{i,j}\in\Lambda$ for $i,j\in{\bf Z}$ by
\begin{equation}\label{IOTA}
\iota_{i,j}=b_{\bar i-2}-b_{\bar j-2}+k\delta(\bar i\leq\bar j)\in\Lambda,
\end{equation}
where $b_{-1}=0$,
$\bar i$ and $\bar j$ are such that $1\leq \bar i,\bar j\leq l$,
$i\equiv \bar i\bmod l$, $j\equiv \bar j\bmod l$, and
\[
\delta(*)=
\begin{cases}
1&\hbox{\rm if $*$ is true};\\
0&\hbox{\rm otherwise}.
\end{cases}
\]
Note that $\iota_{2,1}=b_0\in\Lambda$ and $\iota_{a,b}\not=0$ in $\Lambda$.
In the computation of $x_i$ in (\ref{XI}), we start {}from{} $\iota_{2,1}$,
and apply $c_i,c_{i-1}$, etc., successively.
We will prove that if $M$ is good we will never get $0$ in this process.

Let us start with
\begin{lem}\label{ACTION}
The actions of $M_a$ and $M_b$ on $\Lambda$ satisfy
\begin{eqnarray*}
M_a(\iota_{i,j})
&=&
\begin{cases}
\iota_{i+1-\delta_{\bar i,1},j+1-\delta_{\bar j,1}}
&\hbox{if }(\bar i,\bar j)\not=(1,l);\\
0&\hbox{if }(\bar i,\bar j)=(1,l),
\end{cases}\\
M_b(\iota_{i,j})&=&\iota_{i+1,j+1}.
\end{eqnarray*}
\end{lem}

\begin{proof}
The proof is only case checking.
\end{proof}

\begin{lem}\label{NONZERO}
Let $M=C_1\cdots C_N$ be a good monomial, and $C_i=B$.
We have $M_{c_1}\circ\cdots\circ M_{c_{i-1}}(\iota_{2,1})\not=0$.
\end{lem}
\begin{proof}
Consider the mapping $a,b:\{1,\ldots,l\}\rightarrow\{1,\ldots,l\}$
given by
\begin{eqnarray*}
a(j)=
\begin{cases}
j+1&\hbox{if }j\not=1,l;\\
1&\hbox{otherwise},
\end{cases}\\
b(j)=
\begin{cases}
j+1&\hbox{if }j\not=l;\\
1&\hbox{if }j=l.
\end{cases}
\end{eqnarray*}
We have $M_a(\iota_{j,j'})=\iota_{a(j),a(j')}$ unless $(j,j')=(1,l)$,
and $M_b(\iota_{j,j'})=\iota_{b(j),b(j')}$. Suppose that $C_i=B$.
The trajectory of $1$
by the successive mappings $c_x$, where $x$ changes its values as
$i-1,i-2$, etc., stays at $1$ exactly at the following values of $x$.
The trajectory starts
at $1$ for $x=i$. Its stays at $1$ while $C_x=A$. It changes to $2$
when $C_x=B$ for the first time after $x=i$. Suppose that this is when
$x=i_1$. Then, the trajectory comes back to $1$ when $x=i_1-l+1$.
Since $M$ is good, we have $C_{i_1-l+1}=B$. It stays at $1$ while $C_x=A$,
and changes to $2$ when $C_x=B$ for the first time after $x=i_1-l+1$.
This pattern repeats. In particular, the trajectory comes to $l$ only if
the next $C_x$ is equal to $B$.

On the trajectory of $\iota_{2,1}$ at $x=i$, $0$ appears only if
$M_a$ applies to $\iota_{1,l}$. However, since the trajectory of $1$ at $x=i$
comes to $l$ only if the next mapping is $M_b$, this never occurs.
\end{proof}

Fix a good monomial $M=C_1\cdots C_N$. We define $C_i$, $i\leq0$,
by setting $C_i=B$. For $i\geq1$ such that $C_i=B$ we define
\begin{eqnarray*}
\gamma_i=c_1\circ\cdots\circ c_{i-1}(2),\\
\gamma'_i=c_1\circ\cdots\circ c_{i-1}(1).
\end{eqnarray*}

 Because of Lemma
\ref{ACTION} and \ref{NONZERO}, we can rewrite (\ref{XI}) into
$x_i=\iota_{\gamma_{i+1},\gamma'_{i+1}}$ for $i$ such that $C_{i+1}=B$.
We define $i_-$ to be the
largest integer satisfying $i_-<i$ and $C_{i_-}=B$. Then we have the recursion
relation for $\gamma_i$ and $\gamma'_i$:
\begin{eqnarray}
\gamma'_i&=&\gamma_{i_-},\label{REC1}\\
\gamma_i&=&
\begin{cases}
i+1&\hbox{\rm if }i\leq l-1;\\
\gamma'_{i-l+1}&\hbox{\rm if }i\geq l.\label{REC2}
\end{cases}
\end{eqnarray}
The second relation makes sense because $C_{i-l+1}=B$, which follows {}from{}
$C_i=B$ and the assumption that $M$ is good.

Therefore, we have
\begin{cor}\label{GAMMA}
If $M$ is a good monomial and $C_{i+1}=B$ then 
\[
x_i=\iota_{\gamma_{i+1},\gamma'_{i+1}},
\]
where $\gamma_{i+1}, \gamma'_{i+1}$ are given by \Ref{REC1}, \Ref{REC2}. 
\end{cor}

We defined the equivalence of monomials by $M\sim MA$.
If $(x_0,\cdots,x_{N-1})$ is the extremal configuration for $M$
then the extremal configuration for $MA$ is $(x_0,\cdots,x_{N-1},0)$.
We define the equivalence of configurations so that
$(x_0,\cdots,x_{N-1})\sim(x_0,\cdots,x_{N-1},0)$.
The extremal configuration corresponding to an equivalence class of monomials
is well-defined.

\begin{prop}
Suppose $M=C_1\cdots C_N$ is a good monomial. If a monomial $M'$ gives
a configuration equivalent to that of $M$, then $M'$ is equivalent to $M$.
\end{prop}
\begin{proof}
{}From{} Lemma \ref{NONZERO}, if $M'$ is good its equivalence class
is uniquely determined {}from{} the corresponding configuration.
Suppose $M'=C'_1\cdots C'_{N'}$ is bad. There exists $i\geq1$ such that
$C_i=A$, $C'_i=B$, and $C_j=C'_j$ for $j<i$. Since $C_i=A$ and the
configurations corresponding to $M,M'$ are equivalent, we have
$c_1\circ\cdots\circ c_{i-1}(\iota_{21})=0$. However, this is impossible
because $C_1\cdots C_{i-1}$ is good, and by the same argument as in Lemma
\ref{NONZERO}, we can show that if $c_j\circ\cdots\circ c_{i-1}(1)=l$ then
$C_{j-1}$ is necessarily $B$.
\end{proof}

In Section \ref{sec:5}, we will show that if $M$ is bad, there exists
another bad monomial $M'$ such that $(M+M')v_\infty=0$. Such a pair of
(equivalence classes of) monomials will be 
called a cancellation pair. 
In particular, we will show (Proposition \ref{CANCON})
that the extremal configurations for a cancellation pair are identical.

\setcounter{section}{2}
\setcounter{equation}{0}


\section{Cancellation pairs}\lb{sec:3}
In this section we give a combinatorial preparation 
which will be used to prove Theorem \ref{thm:2.1}.

\subsection{Summation graph}\lb{subsec:3.1}
Let $M$ be a monomial in $A,B$ of degree $m$. 
It acts on a vector $f[P_1,\cdots,P_l]$ as
\begin{eqnarray*}
M(f[P_1,\cdots,P_l])={S}^m(f)\varphi_M(P_1,\cdots,P_l)
[\tilde{P}_1,\cdots,\tilde{P}_l].
\end{eqnarray*}
Here ${S}$ denotes the $q$-shift \eqref{shift}, 
and $\varphi_M(P_1,\cdots,P_l)$ 
is some scalar function determined {}from{} $[P_1,\cdots,P_l]$. 
The vector part 
$[\tilde{P}_1,\cdots,\tilde{P}_l]$ has the structure 
\begin{eqnarray*}
[q^{r_1}z^{s_1}{S}^m(P_{\sigma_M(1)}),
\cdots,q^{r_l}z^{s_l}{S}^m(P_{\sigma_M(l)})] 
\qquad (r_i,s_i\in\Z),
\end{eqnarray*}
where $\sigma_M:\{1,\cdots,l\}\rightarrow \{1,\cdots,l\}$ 
is some map which is not necessarily bijective.  
We have 
\begin{eqnarray}
&&\sigma_A(1)=1,\quad \sigma_A(i)=i-1\quad (2\le i\le l),
\lb{sigma0}\\
&&\sigma_B(1)=l,\quad \sigma_B(i)=i-1\quad (2\le i\le l).
\lb{sigma1}
\end{eqnarray}
Moreover, 
\begin{eqnarray}
\sigma_{MM'}=\sigma_{M'}\circ\sigma_{M}.
\lb{sigma2}
\end{eqnarray}
Note that the map $\sigma_M$ is uniquely 
specified by \eqref{sigma0}--\eqref{sigma2}.

\begin{figure}[h]
\begin{picture}(280,160)(-150,-80)  
\thicklines

\put(15,50){$[\bullet\bullet\bullet]$}
\put(-38,10){$[\bullet\bullet\circ]$}
\put(15,10){$[\bullet\circ\bullet]$}
\put(80,10){$[\circ\bullet\bullet]$}
\put(-38,-30){$[\bullet\circ\circ]$}
\put(15,-30){$[\circ\circ\bullet]$}
\put(80,-30){$[\circ\bullet\circ]$}

\put(51,55){\circle{15}}
\put(46,55){\vector(0,1){0}}
\put(62,57){$B$}

\put(21 ,47){\vector(-3,-2){38}}
\put(-10,35){$A$}

\put(-11,16){\vector(1,0){26}}
\put(2,20){$B$}
\put(44,14){\vector(1,0){28}}
\put(52,17){$B$}
\put(15,10){\vector(-1,0){26}}
\put(3,1){$A$}

\put(73,9){\line(-5,-2){45}}
\put(27,-9){\vector(-3,1){45}}
\put(62,-3){$A,B$}

\put(-28,5){\vector(0,-1){25}}
\put(-38,-3){$A$}


\put(-11,-28){\vector(1,0){26}}
\put(-2,-24){$B$}
\put(42,-28){\vector(1,0){28}}
\put(48,-24){$A,B$}

\put(73,-31){\line(-5,-2){45}}
\put(27,-49){\vector(-3,1){45}}
\put(52,-48){$A,B$}

\put(-49,-29){\circle{15}}
\put(-42,-28){\vector(0,1){0}}
\put(-66,-52){$A$}
\end{picture}
\caption{Summation graph ($l=3$).}
\label{fig:4.1}
\end{figure}
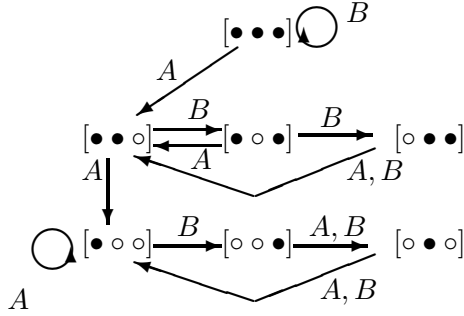

Let us classify monomials according to the subset 
\begin{eqnarray*}
\{\sigma_M(1),\cdots,\sigma_M(l)\}\subset\{1,\cdots,l\}.
\end{eqnarray*}
{}For that purpose it is convenient to introduce an oriented graph $G$, 
which we call the {\em summation graph}.  
By definition, the vertices of $G$ are nonempty subsets of
$\{1,2,\cdots,l\}$. 
We represent a vertex $I$ by an array of length $l$, 
whose $i$-th component is $\bullet$ or $\circ$ 
depending on whether $i\in I$ or not. 
{}For example, $[\bullet\,\circ\,\circ\,\bullet\,\bullet]$ 
means $I=\{1,4,5\}$ with $l=5$. 
Each vertex is a source of two arrows, arrow $A$ and arrow $B$. 
These arrows connect vertices by the following rule. 
\begin{eqnarray}
&&[\bullet\,*\,J]
~~\overset{A}{\longrightarrow}~~
[\bullet\,J\,\circ],
\lb{arrow1}\\
&&[\circ\,J']
~~\overset{A}{\longrightarrow}~~
[J'\,\circ].
\lb{arrow2}\\
&&[*\,J']
~~\overset{B}{\longrightarrow}~~
[J'\,*].
\lb{arrow3}
\end{eqnarray}
Here $*=\bullet,\circ$, and $J,J'$ are 
sequences of $\bullet, \circ$ of length $l-2$ or $l-1$, respectively.
There are only two arrows whose source and sink coincide,
\begin{eqnarray*}
&&[\overbrace{\bullet\,\cdots\,\bullet}^{l}]
~~\overset{B}{\longrightarrow}~~
[\overbrace{\bullet\,\cdots\,\bullet}^{l}],
\\
&&[\bullet\,\overbrace{\circ\,\cdots\,\circ}^{l-1}]
~~\overset{A}{\longrightarrow}~~
[\bullet\,\overbrace{\circ\,\cdots\,\circ}^{l-1}].
\end{eqnarray*}
The summation graph $G$ in the case $l=3$ is shown in fig.\ref{fig:4.1}.

We choose the vertex $I_{top}=[\bullet\,\cdots\,\bullet]$ with all
$\bullet$, 
and fix it as a reference vertex. 
Given a monomial 
$M=C_1C_2\cdots C_m$ ($C_i\in\{A,B\}$), 
let us follow the arrows $C_1,C_2\cdots$ starting {}from{} $I_{top}$. 
We call the resulting sequence of vertices $\{I^{(i)}\}$,  
\begin{eqnarray}
&&I_{top}=I^{(0)}
~~\overset{C_1}{\longrightarrow}~~
I^{(1)}
~~\overset{C_2}{\longrightarrow}~~
\cdots
~~\overset{C_m}{\longrightarrow}~~
I^{(m)},  
\lb{path}
\end{eqnarray}
the {\em path associated with $M$}. 
We say that the path (or, by abuse of language, the monomial $M$)
{\em passes through} the vertices $I^{(i)}$. 
With the rule above, it is easy to see the following. 
\begin{lem}\lb{lem:3.1}
The subset $\{\sigma_M(1),\cdots,\sigma_M(l)\}$
is given by the last vertex $I^{(m)}$ in the path. 
\end{lem}
{}For example, if $l=5$ and $M=A^2B^2A$, then 
\begin{eqnarray}
&&[\bullet\,\bullet\,\bullet\,\bullet\,\bullet]
~~\overset{A}{\longrightarrow}~~
[\bullet\,\bullet\,\bullet\,\bullet\,\circ]
~~\overset{A}{\longrightarrow}~~
[\bullet\,\bullet\,\bullet\,\circ\,\circ]
\lb{ex1}\\
&&\quad 
~~\overset{B}{\longrightarrow}~~
[\bullet\,\bullet\,\circ\,\circ\,\bullet]
~~\overset{B}{\longrightarrow}~~
[\bullet\,\circ\,\circ\,\bullet\,\bullet]
~~\overset{A}{\longrightarrow}~~
[\bullet\,\circ\,\bullet\,\bullet\,\circ],
\nn
\end{eqnarray}
hence $\{\sigma_M(1),\cdots,\sigma_M(5)\}=\{1,3,4\}$.
Notice that in this procedure the monomial should be read {}from{} 
{\em left to right} because of \eqref{sigma2}. 

We say that a vertex $I$ is {\em good} if $1\in I$,  
or equivalently, if the first entry of $I$ is $\bullet$. 
The notion of good monomials \eqref{goodM} 
has the following interpretation. 
\begin{lem}\lb{lem:3.11}
A monomial is good if and only if it passes through 
only good vertices. 
\end{lem}
\begin{proof}
We retain the notation of \eqref{path}. 
If $C_i=A$, then there is a $\circ$ in the $l$-th component 
of $I^{(i)}$.   
This $\circ$ will then move to the left by one step 
with each $A$ or $B$ arrow, 
until it reaches the second component in the array. 
Hence the first component of $I^{(i+l-1)}$ is $\bullet$ if and 
only if $C_{i+l-1}=A$.
\end{proof}

\subsection{Marked vertices}\lb{subsec:3.2}
In dealing with paths associated with monomials, 
it is often convenient to keep track of the move of the last entry of
$I_{top}$. 
Let us introduce paths with marked vertices. 
A marked vertex is a pair $(I,j)$ of vertex $I$ and an element $j\in I$. 
The element $j$ is represented by the symbol $\times$ in place of
$\bullet$.  
We extend the definition of the $A,B$ arrows to marked vertices
as follows.
\begin{eqnarray}
&&[\circ\, J]
~~\overset{A}{\longrightarrow}~~
[J\,\circ],
\lb{arrow4}\\
&&[\bullet\bullet J]
~~\overset{A}{\longrightarrow}~~
[\bullet J\,\circ],
\lb{arrow5}\\
&&[\bullet\circ J]
~~\overset{A}{\longrightarrow}~~
[\bullet J\,\circ],
\lb{arrow6}\\
&&[\bullet\times J]
~~\overset{A}{\longrightarrow}~~
[\bullet J\,\circ],
\lb{arrow7}\\
&&[\times\bullet J]
~~\overset{A}{\longrightarrow}~~
[\bullet J\,\circ],
\lb{arrow8}\\
&&[\times\circ J]
~~\overset{A}{\longrightarrow}~~
[\times J\,\circ],
\lb{arrow9}\\
&&[\bullet\, J]
~~\overset{B}{\longrightarrow}~~
[J\,\bullet],
\lb{arrow10}\\
&&[\circ\, J]
~~\overset{B}{\longrightarrow}~~
[J\,\circ],
\lb{arrow11}\\
&&[\times\, J]
~~\overset{B}{\longrightarrow}~~
[J\,\times].
\lb{arrow12}
\end{eqnarray}
Here $J$ in (\ref{arrow4}--\ref{arrow6}) and (\ref{arrow10}--\ref{arrow11})
contain a $\times$, and others not. Note that in (\ref{arrow7}--\ref{arrow8})
$\times$ is lost in the RHS.

Starting with $[\bullet\,\cdots\,\bullet\,\times]$
we obtain a path $P$ passing through marked or unmarked vertices
similar to \eqref{path}. 
We call it the {\em marked path} associated with $M$.  
{}For example, if $l=5$ and $M=AB^2A$, then 
\begin{eqnarray}
&&
[\bullet\,\bullet\,\bullet\,\bullet\,\times]
~~\overset{A}{\longrightarrow}~~
[\bullet\,\bullet\,\bullet\,\times\,\circ]
~~\overset{B}{\longrightarrow}~~
[\bullet\,\bullet\,\times\,\circ\,\bullet]
\lb{ex2}\\
&&\quad 
~~\overset{B}{\longrightarrow}~~
[\bullet\,\times\,\circ\,\bullet\,\bullet]~~\overset{A}{\longrightarrow}~~
[\bullet\,\circ\,\bullet\,\bullet\,\circ].
\nn
\end{eqnarray}

The rules (\ref{arrow4}--\ref{arrow12}) are so designed that,       
replacing $\times$ with either $\bullet$ or $\circ$,  
we get back to the original rule \eqref{arrow1}--\eqref{arrow3} for 
unmarked vertices. 
The marked path $P$ provides a one-to-one correspondence between 
the unmarked path associated with a monomial $M$ and that of $AM$.  
Indeed, the former is recovered by replacing $\times$ with $\bullet$ in
$P$.
The latter is obtained by first 
supplementing $P$ with a leftmost arrow 
$[\bullet\,\cdots\,\bullet]~~\overset{A}{\longrightarrow}~~$
and then replacing $\times$ with $\circ$. 
For instance, the second procedure applied to 
the example \eqref{ex2} gives \eqref{ex1} for $AM=A^2B^2A$. 
This correspondence will be used extensively in the sequel. 

\subsection{Cancellation pair}\lb{subsec:3.3}
Henceforth we regard an equivalence class of monomials 
as a monomial $M=C_1\cdots C_m$ such that 
$m$ is sufficiently large and $C_i=A$ for sufficiently large $i$. 

Consider the marked path $P$ associated with a good monomial $M$.
If $P$ does not pass through 
a vertex which has $\times$ as the first entry, then $AM$ is good. 
Otherwise, $P$ contains one of the following arrows.
\begin{enumerate}
\item $[\times\,\bullet\,J]~~\overset{A}{\longrightarrow}~~
[\bullet\,J\,\circ]$,
\item $[\times\,\bullet\,J]
~~\overset{B}{\longrightarrow}~~
[\bullet\,J\,\times]$,
\item $[\times\,\circ\,J]
~~\overset{A}{\longrightarrow}~~
[\times\,J\,\circ]$,
\item $[\times\,\circ\,J]
~~\overset{B}{\longrightarrow}~~
[\circ\,J\,\times]$.
\end{enumerate}
The case (iv) does not occur since $M$ is good.  
If (iii) takes place for some arrow $C_j=A$, 
the $J$ must contain $\bullet$. This is because $A$ does not bring
$\times$ to the first column, and when $B$ does so, it brings $\bullet$
back to the last column.
Thus, we can find a $j'>j$ such that $C_{j'}$ is either (i) or (ii). 
Among the arrows of type (i) or (ii), let $C_i$ be such that 
the index $i$ is the largest.  
We call $C_i$ a {\em cancellation $A$-arrow} in the case (i),    
and a {\em cancellation $B$-arrow} in the case (ii). 

Denote by $\mathcal{G}$ the set of equivalence classes 
of good monomials. 
We have shown that it is a disjoint union of three subsets,   
\begin{eqnarray*}
\mathcal{G}=\mathcal{G}_G\cup\mathcal{G}_A\cup\mathcal{G}_B. 
\end{eqnarray*}
Here 
$\mathcal{G}_G=\{M\in\mathcal{G}\mid AM\in\mathcal{G}\}$, 
and $\mathcal{G}_A$ (resp. $\mathcal{G}_B$) is 
the set of good monomials which has a cancellation $A$-arrow (resp.
$B$-arrow). 

\begin{lem}\lb{lem:3.2}
There is a bijection between $\mathcal{G}_A$ and $\mathcal{G}_B$.  
\end{lem}
\begin{proof}
Let $M=C_1\cdots C_m=NAK\in \mathcal{G}_A$ where $N=C_1\cdots C_{i-1}$
and $K=C_{i+1}\cdots C_m$. Suppose that $C_i=A$ is the cancellation $A$-arrow. 
Then, $K$ is a good monomial,   
and the path after $C_i$ has the form 
\begin{eqnarray*}
[\times\,\bullet\,J]
~~\overset{A}{\longrightarrow}~~
[\bullet\,J\,\circ]
~~\overset{C_{i+1}}{\longrightarrow}~~
\cdots
~~\overset{C_{i+l-2}}{\longrightarrow}~~
[\bullet\,\circ\,J']
~~\overset{A}{\longrightarrow}~~
[\bullet\,J'\,\circ]
~~\overset{C_{i+l}}{\longrightarrow}~~
\cdots
\end{eqnarray*}
Replacing the cancellation $A$-arrow by $B$, 
the path associated with $NBK$ after $B$ becomes
\begin{eqnarray*}
[\times\,\bullet\,J]
~~\overset{B}{\longrightarrow}~~
[\bullet\,J\,\times]
~~\overset{C_{i+1}}{\longrightarrow}~~
\cdots
~~\overset{C_{i+l-2}}{\longrightarrow}~~
[\bullet\,\times\,J']
~~\overset{A}{\longrightarrow}~~
[\bullet\,J'\,\circ]
~~\overset{C_{i+l}}{\longrightarrow}~~
\cdots
\end{eqnarray*}
This shows that the $B$ in $NBK$ is the cancellation $B$-arrow. 

Hence we obtain a map {}from{} $\mathcal{G}_A$ into $\mathcal{G}_B$. 
The inverse map can be constructed similarly.
\end{proof}

We call the pair $(NAK,NBK)$ described in the proof 
a {\em cancellation pair}.

\subsection{Reduction of the Theorem}\lb{subsec:3.4}
Let us turn to the proof of Theorem \ref{thm:2.1}.  
Denote by $\chi^*(q,z)$ the right hand side of \eqref{chi*}.
We must show that $\chi^*(q,z)=\chi_{k,l;{\bf b}}(q,z)$. 
Clearly $\chi^*(q,0)=1$.  
Hence, in view of Lemma \ref{infinite recurrence lemma}, 
it is enough to verify the difference equation 
\begin{eqnarray*}
\chi^*(q,z)=(A+B)\chi^*(q,z).
\end{eqnarray*}
The left hand side can be written as 
\begin{eqnarray*}
\sum_{M\in A\mathcal{G}\cap\mathcal{G}}M v_\infty
+\sum_{M\in B\mathcal{G}\cap\mathcal{G}}M v_\infty, 
\end{eqnarray*}
because by the definition if $AM'$ or $BM'$ is a good monomial
then $M'$ is also a good monomial.
Note further that $A\mathcal{G}\cap\mathcal{G}=A\mathcal{G}_G$ and
$B\mathcal{G}\cap\mathcal{G}=B\mathcal{G}$.

The right hand side gives 
\begin{eqnarray*}
\sum_{M'\in\mathcal{G}_G}AM' v_\infty
+\sum_{M'\in\mathcal{G}_A}AM' v_\infty
+\sum_{M'\in\mathcal{G}_B}AM' v_\infty
+\sum_{M'\in\mathcal{G}}BM' v_\infty.
\end{eqnarray*}
Combining these relations with Lemma \ref{lem:3.2}, we see that the
proof of the first statement of 
Theorem \ref{thm:2.1} reduces to the following statement.  
\begin{thm}\lb{thm:3.1}
{}For any cancellation pair $(M,M')=(NAK,NBK)$, we have 
\begin{eqnarray}
AN(A+B)Kv_\infty=0.
\lb{cancel}
\end{eqnarray}
\end{thm}
The proof of Theorem \ref{thm:3.1} will be given in the next section.

Here we give a proof of the following statement.
\begin{prop}\label{CANCON}
Suppose that $(NAK,NBK)$ is a cancellation pair. The extremal configurations
for the pair of monomials $(ANAK,ANBK)$ are identical.
\end{prop}
\begin{proof}
Suppose that $AN=C_1\cdots C_{i-1}$, and $K=C_{i+1}\cdots C_M$.
We are to show that $x_j$ and $x'_j$ 
corresponding to $ANAK$ and $ANBK$, respectively, are in fact equal.
There are four cases (i) $j\leq i-1$, (ii) $j=i-1$, (iii) $j\geq i$ and
$C_{j+1}=A$, and (iv) $j\geq i$ and $C_{j+1}=B$.

The statement is clear if $j\leq i-2$ because the part $AN$ is common
to $ANAK$ and $ANBK$.

For $j=i-1$ since $x_{i-1}=0$,
we must show that $x'_{i-1}=M_{c_1}\circ\cdots M_{c_{i-1}}(\iota_{2,1})=0$.
Consider the marked path $\{I^{(a)}\}$
corresponding to the monomial $ANBK$. We have
$I^{(0)}=[\bullet\cdots\bullet\times]$. Let $1\leq\eta_a,\eta'_a\leq l$ be
such that the position of $\times$ in the marked vertex $I^{(a)}$
is $\eta_a$, and that of the first $\bullet$ we find when we move {}from{}
$\times$ to the right cyclically is $\eta'_a$. It is easy to see that
$\eta_a=c_{a+2}\circ\cdots\circ c_{i-1}(1)$ and
$\eta'_a=c_{a+2}\circ\cdots\circ c_{i-1}(2)$. We have
$(\eta'_{i-2},\eta_{i-2})=(2,1)$ and $(\eta'_0,\eta_0)=(1,l)$.
Therefore, $M_{c_1}\circ\cdots\circ M_{c_{i-1}}(\iota_{2,1})=0$.

Now consider $j\geq i$.
If $C_{j+1}=A$ then $x_j=x'_j=0$.
Take $j\geq i$ such that $C_{j+1}=B$. We will show that
$M_{c_1}\circ\cdots M_a\cdots\circ M_{c_j}(\iota_{2,1})
=M_{c_1}\circ\cdots M_b\cdots\circ M_{c_j}(\iota_{2,1})$,
where $M_a$ in the LHS and $M_b$ 
in the RHS are at the place $i$. Since $NAK$ and $NBK$ are good,
$M_a(\iota_{1,l})=0$ never occurs on the trajectories of $\iota_{2,1}$
in these formulas. Therefore, it is enough to show
$c_1\circ\cdots a\cdots\circ c_j(s)=c_1\circ\cdots b\cdots\circ c_j(s)$
for both $s=1,2$. Note that $a(t)=b(t)$ unless $t=1$. On the other hand,
$c_2\circ\cdots\circ a(1)=l,c_2\circ\cdots\circ b(1)=1$ as shown above.
Therefore, $c_1\circ\cdots\circ a(1)=c_1\circ\cdots\circ b(1)=1$. This proves
the statement.
\end{proof}

\setcounter{section}{3}
\setcounter{equation}{0}


\section{Proof of the cancellation}\lb{sec:4}

\subsection{The pairwise cancellation}\lb{subsec:4.1}
Let us describe the mechanism of the cancellation \eqref{cancel}.
{}First we note a simple fact which is a corollary of the definition of
$A, B$:  
\begin{lem}\lb{lem:4.0}
{}For any simple vector $v$, let 
\begin{eqnarray*}
Av=f[P_1,\cdots,P_l],\quad Bv=f'[P_1',\cdots,P_l'].
\end{eqnarray*}
Then 
\begin{eqnarray*}
&&P_1=P_2,~~P_j=P_j'\quad (2\le j\le l),
\\
&&f/f'=-P_2'/P_1'. 
\end{eqnarray*}
\end{lem}
In view of this Lemma, we set
\begin{eqnarray}
&&BNAKv_\infty=g[Q_1,\cdots,Q_l],
\lb{ve3}\\
&&BNBKv_\infty=g'[Q_1',\cdots,Q_l'],
\lb{ve4}\\
&&AKv_\infty=h[R_1,\cdots,R_l],
\lb{ve1}\\
&&BKv_\infty=h'[R_1',\cdots,R_l'].
\lb{ve2}
\end{eqnarray}
{}{}From{} Lemma \ref{lem:4.0} we have 
\begin{eqnarray*}
&&ANAKv_\infty=-\frac{Q_2}{Q_1}g[Q_2,Q_2,\cdots,Q_l],
\\
&&ANBKv_\infty=-\frac{Q_2'}{Q_1'}g'[Q'_2,Q'_2,\cdots,Q'_l]. 
\end{eqnarray*}

The proof of the following proposition is given in Section \ref{subsec:4.4}.
\begin{prop}\lb{prop:4.1}
Set $m={\rm deg}\,BN$. Then 
\begin{enumerate}
\item $Q_j=Q_j'$ for $2\le j\le l$,
\item $Q_2/Q_1=q^\alpha {S}^m(R_2/R_1)$ and 
$Q'_2/Q'_1=q^\alpha {S}^m(R'_2/R'_1)$, where $\alpha$ is an
integer,
\item $g/g'={S}^m(h/h')$.
\end{enumerate}
\end{prop}
Assuming this proposition we give
\begin{proof} ({\it of Theorem \ref{thm:3.1}})
We have
\begin{eqnarray*}
\left(-\frac{Q_2}{Q_1}g\right)/
\left(-\frac{Q'_2}{Q'_1}g'\right)
={S}^m\left(\frac{R_2}{R_1}\frac{R_1'}{R_2'}\frac{h}{h'}\right)=-1,
\end{eqnarray*}
where in the second equality we used Lemma \ref{lem:4.0}. Therefore
$$
ANAKv_\infty+ANBKv_\infty=0.
$$
\end{proof}

The rest of this section is devoted to the proof of Proposition \ref{prop:4.1}.

\subsection{Cancellation blocks}\lb{subsec:4.2}
In this subsection we study the structure of the monomials $BN$ and $K$ 
when $(NAK,NBK)$ constitutes a cancellation pair. 

Let $s$ be an integer satisfying $0\le s\le l-2$. 
We say that a monomial of degree $l+s$
\begin{eqnarray*}
\mathbb{B}_s=C_1C_2\cdots C_{l+s}
\end{eqnarray*}
is a {\em cancellation block} if the following two conditions
are satisfied.
\begin{enumerate}
\item $C_1=C_{s+2}=C_l=B$,
\item $C_2=\cdots=C_{s+1}=A$, $C_{l+1}=\cdots=C_{l+s}=A$. 
\end{enumerate}

\begin{lem}\lb{lem:4.1}
Let $(M,M')=(NAK,NBK)$ be a cancellation pair. 
Then $BN$ can be decomposed uniquely as a product of cancellation
blocks, 
\begin{eqnarray*}
&&BN=\mathbb{B}_{s_1}\cdots\mathbb{B}_{s_p},
\end{eqnarray*}
with some $p\ge 1$ and $0\le s_1,\cdots,s_p\le l-2$.
\end{lem}
\begin{proof}
First note the following. Suppose a marked vertex has $\times$
at the last column. If we apply only $A$ arrows on this vertex,
we will never reach a cancellation $A$-arrow.

Define $s_1\ge 0$ to be the integer such that $N$ has 
the form $A^{s_1}BN_1$ with some monomial $N_1$. 
The existence of such $B$ in $N$ is assured by the above remark.

We write $BN=C_1C_2\cdots$ with $C_i=A$ or $B$, so that  
$C_1=B$, $C_2=\cdots=C_{s_1+1}=A$ and $C_{s_1+2}=B$.  
Starting with the marked vertex
$[\times\,\bullet\,\cdots\,\bullet]$  
let us keep track of its changes under successive applications 
of $C_1,C_2,\cdots$. 

{}For the first $s_1+1$ letters $BA^{s_1}$ we obtain
\begin{eqnarray*}
[\times\,\bullet\,\cdots\,\bullet]
~~\overset{B}{\longrightarrow}~~
[\bullet\, \cdots\,\bullet\,\times]
~~\overset{A}{\longrightarrow}~~
\cdots
~~\overset{A}{\longrightarrow}~~
[\bullet\, \cdots\,\bullet\,\times\,
\overbrace{\circ\, \cdots \,\circ}^{s_1}].
\end{eqnarray*}

{}From{} the definition (\ref{goodM}) of a good monomial
follows that $s_1\le l-2$ since otherwise $B$ cannot appear in $N$.
Continuing with the move after $BA^{s_1}$
and noting $C_{s_1+2}=B$, we have 
\begin{eqnarray*}
&&[\bullet\, \cdots\,\bullet\,\times\,
\overbrace{\circ\, \cdots \,\circ}^{s_1}]
~~\overset{B}{\longrightarrow}~~
[\bullet\, \cdots\,\bullet\,\times\,
\overbrace{\circ\, \cdots \,\circ}^{s_1}\,\bullet]
\\
&&
~~\overset{C_{s_1+3}}{\longrightarrow}~~
[\bullet\, \cdots\,\bullet\,\times\,
\overbrace{\circ\, \cdots
\,\circ}^{s_1}\,\bullet\,\bar{C}_{s_1+3}]
~~\overset{C_{s_1+4}}{\longrightarrow}~~
\cdots
~~\overset{C_{l-1}}{\longrightarrow}~~
[\bullet\,\times\,
\overbrace{\circ\, \cdots \,\circ}^{s_1}\,\bullet\,\bar{C}_{s_1+3}\,\cdots
\,\bar{C}_{l-1}],
\end{eqnarray*}
where $\bar{C}_i$ signifies $\circ$ if $C_i=A$, and $\bullet$ otherwise. 
The next letter $C_{l}$ must be $B$, since otherwise the $\times$
will disappear by $A$ and it would mean that $ANAK$ is a good monomial.
After that, only $A$ ($s_1$ times) is allowed, and we arrive at
\begin{eqnarray*}
[\bullet\,\times\,
\overbrace{\circ\, \cdots \,\circ}^{s_1}\,\bullet\,\bar{C}_{s_1+3}\,\cdots
\,\bar{C}_{l-1}]
~~\overset{BA^{s_1}}{\longrightarrow}~~
[\times\,\bullet\,\bar{C}_{s_1+3}\,\cdots
\,\bar{C}_{l-1}\,\bullet\, \overbrace{\circ\, \cdots \,\circ}^{s_1}].
\end{eqnarray*}

If the next letter is $A$, this is a cancellation arrow. 
In this case we have shown that $BN=\mathbb{B}_{s_1}$ with 
$\mathbb{B}_{s_1}=BA^{s_1}BC_{s_1+3}\cdots C_{l-1}BA^{s_1}$. 
Otherwise, $BN=\mathbb{B}_{s_1}N'$ with $N'=B\cdots$. 
We can now repeat the same process for $N'$, starting with the vertex
\begin{eqnarray*}
[\times\,\bullet\,\bar{C}_{s_1+3}\,\cdots \,\bar{C}_{l-1}\,\bullet\,
\overbrace{\circ\, \cdots \,\circ}^{s_1}].
\end{eqnarray*}
After a finite number of steps we exhaust $BN$ and arrive at the
cancellation $A$-arrow.
\end{proof}

\begin{lem}\lb{lem:4.2}
In the notation of Lemma \ref{lem:4.1}, 
let $\mathbb{B}_{s_i}=C^{(i)}_1\cdots C^{(i)}_{l+s_i}$. 
Then, for each $1\le i\le p-1$, 
one of the following conditions holds. 
\begin{enumerate}
\item $s_{i+1}=l-2$, 
\item $l\ge s_i+s_{i+1}+3$ and
$C^{(i+1)}_{l-s_i}=C^{(i+1)}_{l-s_i+1}=\cdots
=C^{(i+1)}_{l-1}=A$.
\end{enumerate}
\end{lem}
\begin{proof}
To simplify the notation let us write 
$s=s_i$, $t=s_{i+1}$, 
$\mathbb{B}_s=\mathbb{B}_{s_i}=C_1\cdots C_{l+s}$ 
and $\mathbb{B}_t=\mathbb{B}_{s_{i+1}}=D_1\cdots D_{l+t}$. 
We have seen that $\mathbb{B}_s$ causes a change of vertices 
\begin{eqnarray*} 
[\times\, \bullet\, \cdots]
~~\longrightarrow~~
[\times\,\bullet \bar{C}_{s+3}\, \cdots \,\bar{C}_{l-1} \bullet\, 
\overbrace{\circ\, \cdots \,\circ}^{s}].
\end{eqnarray*} 
We now continue with $\mathbb{B}_t$. 
The leftmost $B$ in $\mathbb{B}_t$ changes the right hand side to 
\begin{eqnarray*} 
&&\phantom{B~}
[\bullet\, \bar{C}_{s+3}\, \cdots\, 
\bar{C}_{l-1}\, \bullet\, \overbrace{\circ\, \cdots \,\circ}^{s}
\,\times].
\end{eqnarray*} 
We distinguish the three cases. 
\begin{enumerate}
\item $t=l-2$,
\item $l\ge s+t+3$,
\item $1\le l-t-2\le s$.
\end{enumerate}
The case (iii) is impossible. 
Indeed, if (iii) takes place, then 
$D_2\cdots D_{t+1}=A^t$ changes the vertex to 
\begin{eqnarray*} 
&&
[\bullet\, \overbrace{\circ\, \cdots\, \circ}^{l-t-2}\,\times\,
 \overbrace{\circ\, \cdots\, \circ}^{t}],
\end{eqnarray*}
and the next letter $D_{t+2}=B$ leads to a bad vertex. 

In the case (ii), $A^t$ followed by 
the next $l-s-t-2$ letters $D_{t+2}\cdots D_{l-s-1}$ change the vertex
to 
\begin{eqnarray*} 
[\bullet\, \overbrace{\circ\, \cdots\, \circ}^{s} \,
\times\,  \overbrace{\circ\, \cdots \,\circ}^{t}\,
\bar{D}_{t+2}\,\cdots \,\bar{D}_{l-s-1}].
\end{eqnarray*}
By the same reason as above, the next $s$ letters
$D_{l-s},\cdots,D_{l-1}$ must all be $A$, the case $s=0$ inclusive.
The proof is over.
\end{proof}

\begin{lem}\lb{lem:4.3}
We retain the notation of Lemma \ref{lem:4.1}.
Then $K$ has the form $\mathbb{D}_1\cdots \mathbb{D}_r$, where 
$\mathbb{D}_i=\mathbb{D}_i' A^{s_p+1}$ and $\mathbb{D}_i'$ has degree
$l-s_p-2$.
\end{lem}
\begin{proof}
Recall the definition (\ref{goodM}).
The proof is straightforward because $NAK$ is a good monomial,
and $N$ contains $A^{s_p}$ at the right end.
\end{proof}

\subsection{Change of vectors}\lb{subsec:4.3}
Our next task is to study how the vectors change under the action of a 
cancellation block.

\begin{lem}\lb{lem:4.4} 
Let $M=C_1\cdots C_l$ be an arbitrary monomial of degree $l$, 
and set 
\begin{eqnarray*}
M\left( f[P_1,\cdots,P_l]\right)
=\tilde{f}[\tilde{P}_1,\cdots,\tilde{P}_l].
\end{eqnarray*}
\begin{enumerate}
\item If $C_i=A$, then $\tilde{P}_i=\tilde{P}_{i+1}$,
\item If $C_i=B$, then $\tilde{P}_i=q^{i-1}z{S}^l(P_i)$.
\end{enumerate}
\end{lem}
\begin{proof} 
{}For each $i=1,\cdots,l$ let 
$C_i\cdots C_l[P_1,\cdots,P_l]=f^{(i)}[P^{(i)}_1,\cdots,P^{(i)}_l]$.
If $C_i=A$, then $P^{(i)}_1=P^{(i)}_2$. 
If $C_i=B$, then
$P^{(i)}_1=z{S}(P^{(i+1)}_l)=z{S}^{l-i+1}(P_i)$. 
Noting that $\tilde{P}_i=P^{(1)}_i={S}^{j-1}(P^{(j)}_{i-j+1})$
($i\ge j$), 
we obtain the assertion.
\end{proof}

Now let $\mathbb{B}_s=C_1\cdots C_{l+s}$ be a cancellation block of
degree $l+s$.  
Given a vector $v=f[P_1,\cdots,P_l]$, we set 
\begin{eqnarray}
\mathbb{B}_s f\left([P_1,\cdots,P_l]\right)
=\tilde{f}[\tilde{P}_1,\cdots,\tilde{P}_l]. 
\lb{ftil}
\end{eqnarray}
\begin{lem}\lb{lem:4.5} 
In the notation of \eqref{ftil}, the following hold 
for $1\le i\le l$.  
\begin{enumerate}
\item If $C_i=A$, then $\tilde{P}_i=\tilde{P}_{i+1}$,
\item If $C_i=B$, then 
\begin{eqnarray*}
\tilde{P}_i=\begin{cases} 
q^{i-1}z{S}^{l+s}(P_1) & (1\le i\le s+1),\\
q^{i-1}z{S}^{l+s}(P_{i-s}) & (s+2\le i\le l).\\
\end{cases}
\end{eqnarray*}
\item $\tilde{P}_i/\tilde{P}_1=q^{s+1}{S}^{l+s}(P_2/P_1)$ ($2\le
i\le s+2$).
\end{enumerate}
\end{lem}
\begin{proof} 
The assertions (i),(ii) follow {}from{} Lemma \ref{lem:4.4} applied to 
\begin{eqnarray*}
A^s\left(f[P_1,\cdots,P_l]\right)
=f'[\overbrace{{S}^{s}(P_1),\cdots,{S}^{s}(P_1)}^{s+1},
{S}^{s}(P_2),\cdots,{S}^{s}(P_{l-s})].
\end{eqnarray*}
Since $C_1=C_{s+2}=B$, we have $\tilde{P}_1=z{S}^{l+s}(P_1)$ and 
$\tilde{P}_{s+2}=q^{s+1}z{S}^{l+s}(P_2)$. 
On the other hand, $C_2=\cdots=C_{s+1}=A$ implies
$\tilde{P}_2=\cdots=\tilde{P}_{s+2}$,  
whence follows the assertion (iii).
\end{proof}

Let us study the scalar part $\tilde{f}$ in \eqref{ftil}.
Define $v^{(j)}$ inductively by $v^{(l+s+1)}=v$ and 
\begin{eqnarray*}
v^{(j)}=C_jv^{(j+1)}=f^{(j)}[P^{(j)}_1,\cdots,P^{(j)}_l].
\end{eqnarray*}
We have $\tilde{f}=f^{(1)}$. 
In general, $f^{(j)}$ has the form 
\begin{eqnarray*}
f^{(j)}=\frac{1}{1-g_j}{S}\left(\frac{1}{1-g_{j+1}}\right)
\cdots
{S}^{l+s-j}\left(\frac{1}{1-g_{l+s}}\right).
\end{eqnarray*}
The factors $g_j$ are determined  as follows. 
\begin{lem}\lb{lem:4.8}
In the above setting, we have
\begin{eqnarray}
&&g_1=q^{s+1}{S}^{l+s}(P_2/P_1),\quad
g_j=q^{j-s-2}{S}^{l+s+1-j}(P_1/P_2)~~(2\le j\le s+1),
\lb{gj1}\\
&&
g_l=z^{-1}{S}^{s+1}(P_1/P_{l-s}),\quad
g_j=z{S}^{l+s+1-j}(P_{j-s}/P_1)~~(l+1\le j\le l+s),
\lb{gj2}\\
&&\mbox{$g_{s+2},\cdots,g_{l-1}$ do not depend on $P_1$}.
\lb{gj3}
\end{eqnarray}
\end{lem}
\begin{proof}
Let us calculate $g_j$'s by applying $C_{l+s},C_{l+s-1},\cdots$ step by
step. 
The formula \eqref{gj2} is easily obtained.
After $s+1$ steps, $[P^{(l-1)}_1,\cdots,P^{(l-1)}_l]$ is given by
\begin{eqnarray*}
[z{S}^{s+1}(P_{l-s}),\overbrace{{S}^{s+1}(P_1),\cdots,{S}^{s+1}(P_1)}^{s+1},
{S}^{s+1}(P_2),\cdots,{S}^{s+1}(P_{l-s-1})].
\end{eqnarray*}
Hence during the next $l-s-2$ steps the computation of $g_j$
involves only 
$P_2,\cdots,P_{l-s}$. 
Since $C_{s+2}=B$, $[P^{(s+2)}_1,\cdots,P^{(s+2)}_l]$ takes the form
\begin{eqnarray*}
[z{S}^{l-1}(P_2),*,\cdots,*,\overbrace{{S}^{l-1}(P_1),\cdots,{S}^{l-1}(P_1)}^{s+1}].
\end{eqnarray*}
The $g_j$ for $1\le j\le s+1$ can be calculated {}from{} this. 
\end{proof}

We say that a pair of simple vectors $(v,v')$ satisfies the 
{\em condition ${\bf (C)}_s$} ($0\le s\le l-2$)
if 
\begin{description}
\item[(C1)$_s$] $P_j=P_j'$  ($2\le j\le l$),
\item[(C2)$_s$] $P_1/P_2=zP_l'/P_1'$,
\item[(C3)$_s$] 
$P_{l-s}=P_{l-s+1}=\cdots=P_l$, 
\end{description}
where $v=f[P_1,\cdots,P_l]$, $v'=f'[P_1',\cdots,P_l']$. 
Because of ${\bf (C1)}_s$,  
these conditions are symmetric in $\{P_j\}_{j=1}^l$
and $\{P'_j\}_{j=1}^l$.

Let $(v, v')$ be a pair satisfying 
${\bf (C)}_s$. 
Let $\mathbb{B}_s$ be a cancellation block, and set  
\begin{eqnarray*}
&&v=f[P_1,\cdots,P_l], \\
&&v'=f'[P_1',\cdots,P_l'],  \\
&&\mathbb{B}_s v=\tilde{f}[\tilde{P}_1,\cdots,\tilde{P}_l], 
\\
&&\mathbb{B}_s v'=\tilde{f'}[\tilde{P}'_1,\cdots,\tilde{P}'_l].
\end{eqnarray*}

\begin{lem}\lb{lem:4.9}
We have
\begin{eqnarray*}
\frac{\tilde{f}}{\tilde{f'}}={S}^{l+s}\left(\frac{f}{f'}\right).
\end{eqnarray*}
\end{lem}
\begin{proof}
As in Lemma \ref{lem:4.8}, define the factors $g_j,g_j'$ 
corresponding to $v,v'$, respectively.
By \eqref{gj2} and {\bf (C3)}$_s$ we have 
\begin{eqnarray}
&&g_l=z^{-1}{S}^{s+1}(P_1/P_l), 
\lb{BBB3}\\
&&g_{l+j}=z{S}^{s+1-j}(P_l/P_1)\quad (1\le j\le s), 
\lb{BBB4}
\end{eqnarray}
and likewise for $g_j'$. 

{}{}From{} \eqref{gj1},\eqref{BBB3} and {\bf (C2)}$_s$ we obtain
\begin{eqnarray*}
{S}^{l-1}(g_l')
=q^{-l+1}z^{-1}{S}^{l+s}(P_1'/P_l')
=q^{s+1}{S}^{l+s}(P_2/P_1)
=g_1.
\end{eqnarray*}
Similarly we find
\begin{eqnarray*}
{S}^{j-1}(g_j)=\begin{cases}
{S}^{j+l-2}(g_{j+l-1}')& (1\le j\le s+1),\\
{S}^{j-l}(g'_{j-l+1}) & (l\le j\le l+s).\\
\end{cases}
\end{eqnarray*}
It is clear {}from{} {\bf (C1)}$_s$ and \eqref{gj3} that 
\begin{eqnarray*}
g_j=g_j'\qquad (s+2\le j\le l-1).
\end{eqnarray*}
The lemma follows {}from{} these relations.
\end{proof}

\begin{lem}\lb{lem:4.7}
Consider a product $\mathbb{B}_s\mathbb{B}_t$ of cancellation blocks 
$\mathbb{B}_s=C_1\cdots C_{l+s}$ and $\mathbb{B}_t=D_1\cdots D_{l+t}$. 
We assume that 
either {\rm (i)} $t=l-2$, 
or {\rm (ii)} $l\ge s+t+3$ and $D_{l-s}=\cdots=D_{l-1}=A$. 
If $(v,v')$ is a pair of vectors satisfying 
${\bf (C)}_t$, 
then the pair $(\tilde{v},\tilde{v}')=(\mathbb{B}_t v,\mathbb{B}_t v')$
satisfies 
${\bf (C)}_s$. 
\end{lem}
\begin{proof}
Set $v=f[P_1,\cdots,P_l]$, $v'=f'[P'_1,\cdots,P'_l]$.  
We check 
${\bf (C)}_s$ for 
$\tilde{v}=\tilde{f}[\tilde{P}_1,\cdots,\tilde{P}_l]$ and 
$\tilde{v}'=\tilde{f'}[\tilde{P}'_1,\cdots,\tilde{P}'_l]$, 
by using Lemma \ref{lem:4.5}.  

{}{}From{} {\bf (C3)}$_s$, Lemma \ref{lem:4.5} (ii) and (iii), we have
\begin{eqnarray}
&&\tilde{P}_1=z{S}^{l+t}(P_1), 
\lb{BBB1}
\\
&&\tilde{P}_l=q^{l-1}z{S}^{l+t}(P_{l-t})=q^{l-1}z{S}^{l+t}(P_{l}),
\lb{BBB2}
\\
&&
\tilde{P}_j/\tilde{P}_1=q^{t+1}{S}^{l+t}(P_2/P_1)\qquad (2\le
j\le t+2).
\lb{AAA}
\end{eqnarray}
The same relations hold by replacing $P_j,\tilde{P}_j$ with
$P_j',\tilde{P}_j'$. 
Combining them with $P_1/P_2=zP_l'/P_1'$ we obtain 
$\tilde{P}_1/\tilde{P}_2=z\tilde{P}_l'/\tilde{P}_1'$,
which proves {\bf (C2)}$_s$. 

Let us show {\bf (C1)}$_s$. 
Note that $D_j=A$ for $2\le j\le t+1$. 
If $D_j=B$ with $j\ge t+2$, then 
$\tilde{P}_j=q^{j-1}z{S}^{l+t}(P_{j-t})$ by Lemma \ref{lem:4.5}
(ii), 
 and 
likewise for $\tilde{P}_j'$.   
Hence {\bf (C1)}$_t$ implies $\tilde{P}_j=\tilde{P}'_j$. 
If $D_i=A$, then we can find some $j$ with $j\ge i+1,t+2$  
such that 
$D_i=D_{i+1}=\cdots=D_{j-1}=A$ and $D_j=B$. 
Therefore, by Lemma \ref{lem:4.5} (i), 
\begin{eqnarray*}
\tilde{P}_i=\tilde{P}_j=\tilde{P}_j'=\tilde{P}_i'.
\end{eqnarray*}

{}Finally, if $t=l-2$, one can show {\bf (C3)}$_s$ directly since
$\mathbb{B}_t=BA^{l-2}B^{l-2}$, and otherwise by using
Lemma \ref{lem:4.5} (i) and the assumption $D_{l-s}=\cdots=D_{l-1}=A$. 
\end{proof}

\subsection{Proof of Proposition \ref{prop:4.1}}\lb{subsec:4.4}
We now prove Proposition \ref{prop:4.1}. 
{}For a cancellation pair $(ANK,BNK)$, 
we use the symbols $g,g',h,h'$ and 
$Q_i,Q_i',R_i,R_i'$ as given in \eqref{ve3}--\eqref{ve2}.
\begin{lem}\lb{lem:4.6}
The vectors $[R_1,\cdots, R_l]$ and 
$[R'_1,\cdots,R'_l]$ satisfy 
${\bf (C)}_{s_p}$.
\end{lem}
\begin{proof}
{}For any $v$ and $D_i$, 
$D_1\cdots D_{l-t-1}A^t v$ has the form $f[P_1,\cdots,P_l]$ with 
$P_{l-t}=\cdots=P_l$. 
Lemma \ref{lem:4.3} implies that 
$Kv_\infty$ has the same form with $t=s_p+1$. 
{}{}From{} this follow {\bf (C2)}$_s$ and {\bf (C3)}$_s$.
The relation {\bf (C1)}$_s$ is clear {}from{} the definition of $A$. 
\end{proof}
\medskip

\noindent{\it Proof of Proposition \ref{prop:4.1}.}\quad
We set $BN=\mathbb{B}_{s_1}\cdots\mathbb{B}_{s_p}$. 
Lemma \ref{lem:4.6} shows that $(v,v')=(AKv_\infty,BKv_\infty)$ 
satisfies ${\bf (C)}_{s_p}$. 
By Lemma \ref{lem:4.7}, we see inductively that the pair
$$
(\tilde{v},\tilde{v}')=
(\mathbb{B}_{s_{i+1}}\cdots\mathbb{B}_{s_p}v, 
\mathbb{B}_{s_{i+1}}\cdots\mathbb{B}_{s_p}v')
$$
satisfies ${\bf (C)}_{s_i}$ for all $i$. 

Statement (i) is an immediate consequence of Lemma \ref{lem:4.7}.
Statement (ii) follows {}from{} Lemma \ref{lem:4.5} (iii) 
by taking $\alpha=\sum_{i=1}^p(s_i+1)$.
{}Finally, {}from{} Lemma \ref{lem:4.9} we conclude that 
$g/g'={S}^m(h/h')$ where $m=\sum_{i=1}^p(s_i+l)$.
This completes the proof of Proposition \ref{prop:4.1}.
\qed

\setcounter{section}{4}
\setcounter{equation}{0}

\def\n{{\bf n}}

\section{Bosonic formulas}\label{sec:5}
In this section we will prove that the sum of contributions
from good monomials is a finite sum of $(l-1)$--fold series
(Corollary \ref{RESULT}). We parameterize the set of good monomials
by a permutation $\sigma$ and a set of integers $(n_1,\ldots,n_{l-1})$
with some restrictions. The vector and scalar part of the corresponding
simple vector is explicitly given in terms of these data.
\subsection{Parameterization of good monomials}
Consider the set $\mathcal{R}$ consisting of
$(\sigma,\n)$ where $\sigma\in\mathfrak{S}_l$
is such that $\sigma(1)=l$ and $\sigma(l)=1$,
and $\n=(n_1,\cdots,n_{l-1})\in{\bf Z}^{l-1}$
satisfying the following conditions:
\begin{eqnarray}
&&n_i\geq0\quad(1\leq i\leq l-2),\label{C1}\\
&&n_{l-1}\geq2-\sigma(l-1),\label{C2}\\
&&n_i>0\quad{\rm if}\quad \sigma(i)<\sigma(i+1).\label{C3}
\end{eqnarray}
\begin{prop}\label{BIJ}
There is a one-to-one correspondence between the set $\mathcal{R}$
and the set of good monomials $\mathcal{G}$.
\end{prop}
\begin{proof}
We choose the representatives of good monomials of the form
$M=C_1\cdots C_n$ where $C_n=B$. First we construct a mapping
$\mathcal{R}\rightarrow\mathcal{G}$. Suppose that $(\sigma,\n)\in\mathcal{R}$.
We define for $i=1,\ldots,l-2$
\begin{equation}
\mathbb{E}_i=C^{(i)}_2\cdots C^{(i)}_l
\end{equation}
where
\begin{equation}
C^{(i)}_j=
\begin{cases}
B&\hbox{\rm if and only if }j\in\{\sigma(1),\cdots,\sigma(i)\};\\
A&\hbox{\rm otherwise}.
\end{cases}
\label{SIGMA}
\end{equation}
We set
\begin{equation}
M=B^{n_{l-1}}\mathbb{E}^{n_{l-2}}_{l-2}\cdots\label{M}
\mathbb{E}^{n_1}_1.
\end{equation}

When $n_{l-1}<0$ we understand this monomial in the following sense.
Set
\begin{equation}
r={\rm max}\{a;1\leq a\leq l-2,n_a\not=0\}.
\label{MAXR}
\end{equation}
Because of (\ref{C3}) we have $\sigma(r+1)>\cdots>\sigma(l-1)$
and $C^{(r)}_j=A$ if and only if $j$ is equal to one of these $l-r-1$
integers. We have, in particular, that the first (i.e., left)
$\sigma(l-1)-2$ elements in $\mathbb{E}_r$ are $B$. Therefore, by (\ref{C2})
we can find
$M$ such that $B^{-n_{l-1}}M=\mathbb{E}_r^{n_r}\cdots\mathbb{E}_1^{n_1}$.

In order to show that $M$ given by (\ref{M}) is good
it is enough to show that if $C^{(i-1)}_j=B$ then $C^{(i)}_j=B$.
This is clear because the only place where $C^{(i-1)}_j\not=C^{(i)}_j$
occurs is at $j=\sigma(i)$, i.e., we have $C^{(i-1)}_{\sigma(i)}=A$
and $C^{(i)}_{\sigma(i)}=B$.

Next, we construct the inverse mapping
$\mathcal{G}\rightarrow\mathcal{R}$.

Suppose $M=C_1\cdots C_n\in\mathcal{G}$. We modify $M$ to $M'=B^{l-1}M$, which
is still good. We construct $(\sigma,\n)$ for $M'$ by the procedure given below.
Then, the one corresponding to $M$ is obtained
{}from{} $(\sigma,\n)$ by changing $n_{l-1}$ to $n_{l-1}-l+1$.

We write $M'=B^{m}\mathbb{D}_t\cdots\mathbb{D}_1$
where $m\geq0$ and $\mathbb{D}_s=D^{(s)}_2\cdots D^{(s)}_l$ and
all $\mathbb{D}_s$ contains at least one $A$. This is, in general, possible
because of the modification of $M$ to $M'$ in the first step.

Because of the condition (\ref{goodM}), we have $D^{(s)}_j=B$
if $D^{(s-1)}_j=B$. For convenience we define $\mathbb{D}_{t+1}=B^{l-1}$.
Namely, $D^{(t+1)}_j=B$ for all $2\leq j\leq l$. 
We say $D^{(s)}_j=B$ is new if
$D^{(s-1)}_j=A$. Besides, all $B$ in $\mathbb{D}_1$ are new by definition.
We define $\sigma(i)$ for $1\leq i\leq l-1$ so that
new $B$'s in $\mathbb{D}_{t+1}\cdots\mathbb{D}_1$,
when they are read {}from{} right to left, are
$D^{(s_1)}_{\sigma(1)},\cdots,D^{(s_{l-1})}_{\sigma(l-1)}$
for some $s_i$ such that $s_1\leq\cdots\leq s_{l-1}$.
We set $\sigma(l)=1$. Note that $\sigma(1)=l$ because $C_n=B$.
Thus we obtained an element $\sigma$ of $\mathfrak{S}_l$.

Now, $\mathbb{E}_i$ and $n_i$ for $1\leq i\leq l-2$
are uniquely determined by the condition that
$\mathbb{D}_t\cdots\mathbb{D}_1$ is of the form
$\mathbb{E}^{n_{l-2}}_{l-2}\cdots\mathbb{E}^{n_1}_1$
where (\ref{SIGMA}) is satisfied.

We have $M'=B^{m}\mathbb{E}^{n_{l-2}}_{l-2}\cdots\mathbb{E}^{n_1}_1$.
Therefore, we have
$M=B^{m-l+1}\mathbb{E}^{n_{l-2}}_{l-2}\cdots\mathbb{E}^{n_1}_1$,
and define $n_{l-1}=m-l+1$.

By the construction the conditions (\ref{C1}) and (\ref{C3})
are satisfied. Note that the first $\sigma(l-1)-2$ elements of
$\mathbb{D}_t\cdots\mathbb{D}_1$ are $B$.
Since $M'=B^{l-1}M$, it implies $\sigma(l-1)-2+m\geq l-1$.
Therefore, the condition (\ref{C2}) is satisfied. 

We have constructed mappings $\mathcal{R}\rightarrow\mathcal{G}$
and $\mathcal{G}\rightarrow\mathcal{R}$. 
By comparing two constructions we can prove that
they are inverse to each other.
\end{proof}

\def\S{{S}}

\subsection{Actions of good monomials on the vector part}
The aim of this section is to describe the actions of good monomials
$M$ given by (\ref{M}) where we consider $A$ and $B$ as operators
acting on $l$-component vectors $[P_1,\cdots,P_l]$ by
\begin{eqnarray*}
A:[P_1,\cdots,P_l]&\mapsto&[{S}(P_1),{S}(P_1),{S}(P_2),
\cdots,{S}(P_{l-1})],\\
B:[P_1,\cdots,P_l]&\mapsto&[z{S}(P_l),{S}(P_1),\cdots,
{S}(P_{l-1})].
\end{eqnarray*}
Namely, these actions are are obtained {}from{} that of $A$ and $B$
forgetting the scalar part. We call it the action on the vector part
and denote by the same letters $A,B$.

First we describe the action of $B^m$ on the vector part.
Suppose $[P_1,\cdots,P_l]$ is given. Define $P^{(n)}_i$
for all $i,n\in{\bf Z}$ by
\begin{eqnarray}
P_{j+ml}&=&q^{-\frac12m(m-1)l-(j-1)m}z^{-m}P_j,\label{P1}\\
P^{(n)}_i&=&\S^n(P_{i-n}).\label{P2}
\end{eqnarray}

\begin{lem}\label{LEMA}
The formula \Ref{P1} is consistent.
The action of $B^n$ $(n\in{\bf Z})$ on the vector part is given by
\begin{equation}
B^n[P_1,\cdots,P_l]=[P^{(n)}_1,\cdots,P^{(n)}_l].
\label{BN}
\end{equation}
\end{lem}
\begin{proof}
We define $P_i$ for arbitrary $i=j+ml\in{\bf Z}$ by (\ref{P1})
where $1\leq j\leq l$ and $m\in{\bf Z}$. The consistency of (\ref{P1})
can be checked by induction on $j$ for all $j$, reducing the case $j>l$ to 
$j-l$, and $j<0$ to $j+l$. We then define $P^{(n)}_i$ 
for all $i,n\in{\bf Z}$ by (\ref{P2}).
The validity of (\ref{BN}) can be checked by induction on $n$,
reducing $n>0$ to $n-1$, and $n<0$ to $n+1$.
\end{proof}

The following formula follows {}from{} $\S\circ z=qz\circ\S$.
\begin{equation}
P^{(n)}_i=q^{\frac12m(m+1)l+(1-i)m}z^{-m}\S^n(P_j)
\quad{where}\quad i-n=j+ml.
\end{equation}

Next we describe the action of $\mathbb{E}_i^{n_i}$ on the vector part.
For this purpose we prepare a few notations. We consider a subset
$I\subset[2,l]$ such that $l\in I$. Let $r$ be the cardinality of $I$.
We identify the subset $I$ with a mapping $I:{\bf Z}\rightarrow[1,l]$ by
\begin{eqnarray*}
&&I(0)=1,\\
&&I=\{I(1),\cdots,I(r)\}\quad{where}\quad I(1)<\cdots<I(r),\\
&&I(i+r+1)=I(i).
\end{eqnarray*}
Note that $I(r)=l$. If $r=l-1$, we have $I(i)=i+1$ for $0\leq i\leq l-1$.

We also define $I^*:[1,l]\rightarrow[0,r]$ by
\begin{equation}\label{INV}
I^*(j)={\rm max}\{i;0\leq i\leq r,I(i)\leq j\}.
\end{equation}
Note that $I^*\circ I(i)=i$ $(0\leq i\leq r)$, $I^*(1)=0$ and
$I^*(l)=r$. If $r=l-1$, we have $I^*(j)=j-1$.

Let $I$ be as above. We define an operator $B[I]=C_2\cdots C_l$
where $C_i=A$ or $B$ and $C_i=B$ if and only if $i\in I$.
As we will see in the following proposition, the action of $B[I]^n$
on the vector part $[P_1,\cdots,P_l]$ depends only on $P_i$ $(i\in I)$,
and the resulting vector part $[P'_1,\cdots,P'_l]$ has the property
$P'_j=P'_{I\circ I^*(j)}$. Namely, the action of $B[I]^n$ on the vector part
is actually understood as acting on the $r+1$ component vector
$P_{I(i)}$ ($0\leq i\leq r$). This is the reason for introducing the mappings
$I$ and $I^*$.
\begin{lem}\label{LEMB}
Let $I$ and $r$ be as above and
$[P'_1,\cdots,P'_l]=B[I]^n[P_1,\cdots,P_n]$ $(n\in{\bf Z}_{\geq0})$.
We have
\begin{eqnarray}
P'_j&=&q^{\omega_r(I,n,i)}z^{\rho_r(n,i)}\S^{(l-1)n}(P_{I(i+n)})
\hbox{ where }i=I^*(j),\label{VEC}\\
\omega_r(I,n,i)
&=&\sum_{s=1}^n(I(i+s)-2)+\frac12n(n-1)(l-1)\nonumber\\
&&\quad+m-(l-1)\{(r-i)m+\frac12m(m-1)(r+1)\},\\
\rho_r(n,i)&=&n-m,
\end{eqnarray}
where $m=m(n,i)$ denotes the integer part of $(i+n)/(r+1)$.
\end{lem}
\begin{proof}
For $n=1$ we have
\begin{eqnarray*}
&&P'_{I(i)}=\cdots=P'_{I(i+1)-1}=q^{I(i+1)-2}z\S^{l-1}(P_{I(i+1)})
\quad(0\leq i\leq r-1),\\
&&P'_{I(r)}(=P'_l)=\S^{l-1}(P'_1).
\end{eqnarray*}
Set
\[
S(i)=
\begin{cases}
q^{I(i)-2}z\S^{l-1}&\hbox{ if }i\not\equiv0\bmod r+1;\\
\S^{l-1}&\hbox{ otherwise}.
\end{cases}
\]

The general case reads as
\[
P'_{I(i)}=S(i+1)\circ\cdots\circ S(i+n)(P_{I(i+n)}).
\]
The result (\ref{VEC}) is obtained by using the commutation relations
$\S\circ z=qz\circ\S$. The number of times that $i+s\equiv0\bmod r+1$ occurs
for $1\leq s\leq n$ is given by $m(n,i)$. For such $s$ the counting of
the powers of $q$ and $z$ differs {}from{} $i+s\not\equiv0\bmod r+1$.
Taking this difference into account we obtain (\ref{VEC}).
\end{proof}

\begin{prop}
Let $M$ be a good monomial and $[P_1,\cdots,P_l]$ 
the vector part of $M[1,\cdots,1]$.
Let $(\sigma,\n)$ correspond to $M$ as in Proposition \ref{BIJ}.
Fix $1\leq i\leq l$ and define $j$ and $m$ by
$i-n_{l-1}=j+ml$ and $1\leq j\leq l$. Define
$I_a, i_a$ $(1\leq a\leq l-2)$ by
$$
I_a=\{\sigma(1),\cdots,\sigma(a)\} \qquad
i_a=
\begin{cases}
I_{l-2}^*(j)&\hbox{ if }a=l-2;\\
I_a^*\circ I_{a+1}(i_{a+1}+n_{a+1})&\hbox{ if }1\leq a\leq l-3,
\end{cases}
$$
(see (\ref{INV})). Then, we have
\begin{eqnarray}
&&P_i=q^{\omega(\sigma,\n)_i}z^{\rho(\sigma,\n)_i},\\
&&\omega(\sigma,\n)_i
=\frac12m(m+1)l+(1-i)m+\sum_{a=1}^{l-2}\omega_a(I_a,n_a,i_a)\nonumber\\
&&\quad+n_{l-1}\sum_{a=1}^{l-2}\rho_a(n_a,i_a)+(l-1)\sum_{1\leq a<b\leq l-2}
n_b\rho_a(n_a,i_a),\\
&&\rho(\sigma,\n)_i=-m+\sum_{a=1}^{l-2}\rho_a(n_a,i_a).
\end{eqnarray}
\end{prop}
\begin{proof}
This proposition follows from Lemma \ref{LEMA} and Lemma \ref{LEMB}.
\end{proof}

Consider $\omega(\sigma,\n)_i$ and $\rho(\sigma,\n)_i$ as polynomials of
$\n=(n_1,\ldots,n_{l-1})$.
They depend on $\sigma$ and $i$. However,
the quadratic part $\omega^{(2)}(\n)$ in $\omega(\sigma,\n)_i$ and
the linear part $\rho^{(1)}(\n)$ in $\rho(\sigma,\n)_i$ are independent:
\begin{eqnarray*}
&&\omega^{(2)}(\n)=\frac{l-1}2\sum_{1\leq a,b\leq l-2\atop c={\rm min}(a,b)}
n_an_b\frac c{c+1}
+n_{l-1}\sum_{1\leq a\leq l-2}n_a\frac a{a+1}+\frac{n_{l-1}^2}{2l},\\
&&\rho^{(1)}(\n)=
\sum_{1\leq a\leq l-2}n_a\frac a{a+1}+\frac{n_{l-1}}l.
\end{eqnarray*}

\subsection{Coloring blocks and scalar factors}
In this section, we determine the scalar part of
the monomial $Mv_\infty$. It is of the form
\[
\frac{(-1)^\alpha q^\beta z^\gamma}
{(q)_{t_2}\cdots(q)_{t_l}(q^{\sum_{c=2}^lt_c}z)_\infty}
\]
We will determine $\alpha,\beta,\gamma$ and $t_a$ ($2\leq l$)
for each $M=C_1\cdots C_n$.

A sub-monomial of the form $C_j\cdots C_{j+r-1}$
($1\leq r\leq l-1,1\leq j\leq n+1-r$) is called a block if 
$j=n,r=1$ or the following conditions are satisfied:
\[
C_i=
\begin{cases}
B&\hbox{for $i=j,j+r$};\\
A&\hbox{for $j+1\leq i\leq j+r-1$.}
\end{cases}
\]
Namely, a block is a maximal (i.e., $C_{j+r}=B$) submonomial
of the form $BA\cdots A$. The restriction $r\leq l-1$ is a consequence
of the assumption that $M$ is good and $C_n=B$.

Set ${\bf B}=\{j;1\leq j\leq n,C_j=B\}$.
Obviously, the blocks in $M$ are labeled by ${\bf B}$.
We define the length and color of each block by mappings
$r:{\bf B}\rightarrow\{1,\ldots,l-1\}$ and
$c:{\bf B}\rightarrow\{1,\ldots,l\}$. The length $r(j)$ of
the block $C_j\cdots C_{j+r-1}$ is $r$. The color of each block
is defined inductively as follows. We set $c(n)=1$. Suppose
we have defined $c(j)$ for all $j\in{\bf B}$ such that $i<j$.
Set $i_1=i+r(i)\in{\bf B}$. If $i_1+l-1>n$ or $C_{i_1+l-1}=A$,
we set $c(i)={\rm max}\{c(j);j\in{\bf B},i<j\}+1$. Otherwise,
we have $i_1+l-1\in{\bf B}$; we set $c(i)=c(i_1+l-1)$.

It is convenient to show the color $c$ of $B$ in $M$ by $B_c$.
For example, if $l=4$, the coloring of $M=B^5(BAB)^4(AAB)^3$
is
\[
B_3B_1B_4B_2B_3B_1AB_2B_3AB_1B_2AB_3B_1AB_2AAB_1AAB_2AAB_1.
\]

Fix a monomial $M$. Define $\lambda:{\bf B}\rightarrow{\bf Z}_{\geq0}$
as follows.
\begin{eqnarray}
\lambda_i=
\begin{cases}
\displaystyle
\sum_{j\in{\bf B}\atop i<j,c(j)\not=1}r(j)&\hbox{if $c(i)=1$};\\
\displaystyle
\sum_{j\in{\bf B}\atop i\leq j,c(j)=c(i)}r(j)&\hbox{if $c(i)\not=1$}.
\end{cases}
\end{eqnarray}
Namely, if $c(i)=1$, $\lambda_i$ is the total length of
the blocks of color not equal to $1$ that are located to the right of
the block $i$. If $c(i)\not=1$, $\lambda_i$ is the total length
of the blocks of color $c(i)$ that are located to the right of
the block $i$ including itself.

\begin{lem}\label{NARROW}
If $C_i=C_j=B_c$ where $i<j$, then we have $i+l\leq j$.
\end{lem}
\begin{proof}
We prove the statement by induction on $i$ starting {}from{} $i=n$
where the statement is valid because the ``if'' statement is false.
Suppose $C_i=C_j=B_c$ where $i<j$ and $i+l>j$. Let $i'$ be the smallest integer
such that $i<i'$ and $C_{i'}=B$. Since the color $c$ is attained at $j$
before $i$, by the way of coloring we have
$i'+l-1\leq n$ and $C_{i'+l-1}=B_c$. Note that $i'\leq j<i+l\leq i'+l-1$,
and therefore we have $0<i'+l-1-j\leq l-1$ and $C_j=C_{i'+l-1}=B_c$,
which contradicts to the induction hypothesis.
\end{proof}

The following proposition determines the scalar factors.

\begin{prop}
Let $M$ be a good monomial, $M=C_1\cdots C_n$ and $C_n=B$.
Set $M_i=C_i\cdots C_n$ and $M_iv_\infty=f_i[P^{(i)}_1,\ldots,P^{(i)}_n]$.
We set formally $M_{n+1}v_\infty=v_\infty=f_{n+1}[1,\ldots,1]$,
where $f_{n+1}=\frac1{(z)_\infty}$. For $1\leq i\leq n$, define
$i_b\in{\bf B}$ such that $i_b\leq i\leq i_b+r(i_b)-1$.
Then, we have
\begin{equation}\label{FAC}
h_i=\frac{f_i}{\S(f_{i+1})}=
\begin{cases}
1/(1-q^{\lambda_{i_b}}z)
&\hbox{if }C_i=A\hbox{ and }c(i_b)=1;\\
1/(1-q^{-\lambda_{i_b}}z^{-1})
&\hbox{if }C_i=B\hbox{ and }c(i_b)=1;\\
1/(1-q^{-\lambda_{i_b}+i-i_b})
&\hbox{if }C_i=A\hbox{ and }c(i_b)\not=1;\\
1/(1-q^{\lambda_{i_b}-i+i_b})
&\hbox{if }C_i=B\hbox{ and }c(i_b)\not=1.\\
\end{cases}
\end{equation}
\end{prop}
\begin{proof}
We prove the statement by induction on $i_b$ starting {}from{} $i_b=n$.
Our induction hypotheses are the following.

(i) The vector part of $M_{i_b+r(i_b)+1}v_\infty$ is of the form
$[q^{\alpha_1}:\cdots:q^{\alpha_l}]$ if $c(i_b)=1$,
otherwise it is $[q^{\alpha_1}z:\cdots:q^{\alpha_l}]$;

In the former case, all components
are of the form $q^{\alpha_j}$, and in the latter there exists $1\leq j_0<l$
such that $q^{\alpha_j}z$ if $j\leq j_0$ and $q^{\alpha_j}$ if $j>j_0$.

(ii) $\alpha_l-\alpha_1=\begin{cases}\lambda_i&\hbox{if }c(i_b)=1;\\
-\lambda_i&\hbox{otherwise}.\end{cases}$.

In the course of induction, we will also prove (\ref{FAC}).

If $i_b=n$, then $i=n$ and $\lambda_n=0$. We have
$c(n)=1$ and $h_n=1-z^{-1}$ by (\ref{An}), (\ref{v infinity})
and (\ref{OPB}). Note also that the vector part of $M_{n+1}v_\infty$
is $[1,\cdots,1]$.

In the below we use the projective notation
$[Q_1:\cdots:Q_l]$ for the vector part $[P_1,\ldots,P_l]$. Namely,
$Q_i=QP_i$ where $Q$ is independent of $i$. 
We use the symbols $\uparrow$, $\downarrow$ and $|$ to indicate the
submonomial sitting in the right of these symbols.
For example, if we write $C_1\cdots C_i\uparrow C_{i+1}\cdots C_n$,
the vector part at $\uparrow$ means that of 
$C_{i+1}\cdots C_nv_\infty$.

Consider $i$ such that $c(i_b)=1$. The monomial $M_{i_b}$ is of the form
\[
\underbrace{B_1A\cdots A}_s\uparrow
\underbrace{B\cdots
\overbrace{BA\cdots A}^r}_{l-1}B_1\downarrow C_{i+s+l}\cdots C_n.
\]
Here $A\cdots A$ are all $A$, but $B\cdots B$ can be a mixture.
We have the restriction $s\leq r$ because $M$ is good.
Note that $C_i$ belongs to $\underbrace{B_1A\cdots A}_s$.

By the induction hypothesis the vector part is of the form
$[q^{\alpha_1}:\cdots:q^{\alpha_l}]$ at $\downarrow$.
Then, by the following $B_1$, the vector part changes to
$[q^{\alpha_l}z:q^{\alpha_1}:\cdots]$, and after the next $r-1$ of $A$,
it further changes to
\[
[\underbrace{q^{\alpha_l+r-1}z:\cdots:q^{\alpha_l+r-1}z}_r:
q^{\alpha_1}:\cdots].
\]
At the point $\uparrow$ it is of the form
\[
[q^{\alpha_1}:\cdots:
\underbrace{q^{\alpha_l+l-1}:\cdots:q^{\alpha_l+l-1}}_r].
\]

The assertion (\ref{FAC}) and the statement (i) immediately follow from this.

Note that the difference $\alpha_l-\alpha_1$ have increased
by $l-1$ {}from{} $\downarrow$ to $\uparrow$.
By Lemma \ref{NARROW} there is only one $B_1$ between $\uparrow$ and
$\downarrow$. Therefore, the total length of blocks of color not $1$
between $\uparrow$ and $\downarrow$ is equal to $l-1$. Therefore,
the statement (ii) is also proved.

Next, we consider the case $c(i_b)=c$ where $c\not=1$.
There are two cases: (a) $C_j\not=B_c$ for all $j>i_b$, and
(b) there exists $C_j=B_c$ for some $j>i_b$. In the latter case, we have
$i_b<j-l+1$, $C_{j-l+1}=B$ and $C_{i'}=A$ for $i_b<i'<j-l+1$.

Case (a): the monomial $M_i$ is of the form
\[
\underbrace{B_c\uparrow A\cdots A\overbrace{|B\cdots\downarrow}^t}_s
\overbrace{B\cdots\underbrace{BA\cdots A}_r}^{l-1}B\cdots.
\]
Here $A\cdots A$ are all $A$, but other dots can be $A$ or $B$.
We have $s\leq r$ because $M$ is good. $C_i$ belongs to the block
$\underbrace{B_c\uparrow A\cdots A}_{s-t}$.

At the point $\downarrow$
the vector part is of the form $[\cdots:\underbrace{1:\cdots:1}_r]$.
It changes to $[z:\cdots:\underbrace{1:\cdots:1}_{r-t}]$ at $|$, and then
to $[q^{s-t-1}z:\cdots:\underbrace{1:\cdots:1}_{r-s}]$ at $\uparrow$.
The assertion (\ref{FAC}) and statements (i) and (ii) follow from these
observations.

Case (b): the monomial $M_i$ is of the form
\[
\underbrace{B_c\uparrow A\cdots A|}_s
\overbrace{B\cdots\underbrace{BA\cdots A}_r}^{l-1}B_c\downarrow\cdots,
\]
where $s\leq r$ and there are no $B_c$ between the two $B_c$ in this diagram.
By the induction hypothesis, the vector part at the point $\downarrow$ is
of the form $[q^{\alpha_1}z:\cdots:q^{\alpha_l}]$. It changes to
$[q^{\alpha_1+1}z:\cdots:\underbrace{q^{\alpha_l}:\cdots:q^{\alpha_l}}_r]$
at $|$,
and then to $[\underbrace{q^{\alpha_1+s}z:\cdots:q^{\alpha_1+s}z}_s
:\cdots:\underbrace{q^{\alpha_l}:\cdots:q^{\alpha_l}}_{r-s+1}]$
at $\uparrow$. The assertion (\ref{FAC}) and the statements (i) and (ii)
follow from these observations.
\end{proof}

For a given good monomial $M=C_1\cdots C_n$ with $C_n=B$ we define $t_c$
and $m_c$ $(1\leq c\leq l)$: the total length of color $c$ blocks is denoted by
$t_c$, and the number of color $c$ blocks is denoted by $m_c$. We also set
\[
\mu_c=\sum_{i:C_i=B_c}(n-i+1).
\]
\begin{cor}\label{RESULT}
The scalar part $f$ of
$Mv_\infty=C_1\cdots C_nv_\infty=f[P_1,\ldots,P_l]$ is given by the formula
\begin{eqnarray}
f&=&\frac{(-1)^\alpha q^{\bar{\beta}} z^{m_1}}
{(q)_{t_2}\cdots(q)_{t_l}(q^{\sum_{c=2}^lt_c}z)_\infty},\\
\alpha&=&n-t_1+\sum_{c=1}^lm_c,\\
\bar\beta&=&nm_1+\frac12\sum_{c=2}^lt_c(t_c+1)-\sum_{c=1}^l\mu_c.
\end{eqnarray}
\end{cor}

As in the case of vector part, the quadratic part of $\bar\beta$
and the linear parts of $m_1$, $t_2,\ldots,t_l$ are independent of $\sigma$.
This follows from the following formula for the linear parts.
Below $x\sim y$ denotes an equality modulo constant term in $n_1,\ldots,n_{l-2}$.
\begin{eqnarray*}
n&=&(l-1)\sum_{a=1}^{l-2}n_a+n_{l-1},\\
m_c&\sim&\sum_{a=c-1}^{l-2}\frac a{a+1}n_a+\frac{n_{l-1}}l,\\
t_c&\sim&(l-1)\sum_{a=c-1}^{l-2}\frac{n_a}{a+1}+\frac{n_{l-1}}l.\\
\end{eqnarray*}

\bigskip

\noindent
{\it Acknowledgments.}\quad 
This work is partially supported by
the Grant-in-Aid for Scientific Research (B)
no.12440039 and (A1) no.13304010, Japan Society for the Promotion of Science.
The work of SL is in part supported by grant RFBR-01-01-00546


\end{document}